\documentclass[12pt,psamsfonts]{amsart}

\newtheorem{thm}{Theorem}

\newtheorem{pro}{Proposition}
\newtheorem{lem}{Lemma}

\newcommand{\rn}{{\mathbb{R}}^n}

\newcommand{\rkn}{{\mathbb{R}}^{n-k}}
\newcommand{\eps}{\varepsilon}
\RequirePackage{amsfonts}

 \RequirePackage{amssymb}
\RequirePackage{latexsym}

\begin{document}
\thispagestyle{empty} \setcounter{page}{1}

\title[ Lipschitz generalization of Malkin-Loud
result] 
{Bifurcations from nondegenerate families of periodic solutions
}

\author[Adriana Buic\u{a}, Jaume Llibre and Oleg Makarenkov]
{Adriana Buic\u{a}, Jaume Llibre and Oleg Makarenkov}

\address{Department of Applied Mathematics,
Babe\c{s}-Bolyai University, Cluj-Napoca, Romania}

\email{abuica@math.ubbcluj.ro}

\address{Departament de Matem\`{a}tiques, Universitat Aut\`{o}noma de
Barcelona, 08193 Bellaterra, Barcelona, Catalonia, Spain}

\email{jllibre@mat.uab.cat}

\address{Department of Mathematics, Voronezh State University, Voronezh, Russia}

\email{omakarenkov@math.vsu.ru}

\thanks{The first   author is partially supported by a MCYT/FEDER grant number MTM2008-00694.
The second author is partially supported by a MEC/FEDER grant number
MTM2008-03437 and by a CICYT grant number 2005SGR 00550. The third
author is partially supported by the Grant BF6M10 of Russian
Federation Ministry of Education and CRDF (US), by the RFBR Grants
09-01-92429, 09-01-00468 and by the President of Russian Federation Young PhD
Student grant MK-1620.2008.1.}


\date{}
\dedicatory{}

\begin{abstract} By a nondegenerate $k$-parameterized family $K$ of periodic solutions we understand the situation
when the geometric multiplicity of the multiplier $+1$ of the
linearized on $K$ system equals to $k.$ Bifurcation of
asymptotically stable periodic solutions from $K$ is well studied
in the literature and different conditions have been proposed
depending on whether the algebraic multiplicity of $+1$ is $k$ or
not (by Malkin, Loud, Melnikov, Yagasaki). In this paper we
assume that the later is unknown. Asymptotic stability can not be
understood in this case, but we demonstrate that the information
about uniqueness of periodic solutions is still available.
Moreover, we show that differentiability of the right hand sides
is not necessary for the results of this kind and our theorems are
proven under a kind of Lipschitz continuity.
\end{abstract}
\maketitle

\section{Introduction}

In  \cite{mal}  Malkin studied the bifurcation of $T$-periodic
solutions in the $n$--dimensional $T$-periodic systems of the form
\begin{equation}\label{mig}
  \dot x = f(t,x)+\varepsilon g(t,x,\varepsilon),
\end{equation}
where both functions $f$ and $g$ are sufficiently smooth.
It is assumed in \cite{mal} that the unperturbed system (namely
(\ref{mig}) with $\varepsilon=0$) has a family of $T$-periodic
solutions, denoted $x(\cdot, \xi(h))$, whose initial conditions are
given by a smooth function $\xi:\mathbb{R}^k\to\mathbb{R}^n$,
$h\mapsto \xi(h)$. In these settings  the adjoint linearized
differential system
\begin{equation}\label{nnn}
  \dot u=-\left(D f(t,x(t,\xi(h)))\right)^*u
\end{equation}
has $k$ linearly independent $T$-periodic solutions $u_1(\cdot,h),
\dots , u_k(\cdot,h)$ and the geometric multiplicity of the
multiplier $+1$ of (\ref{nnn}) is, therefore, $k.$ Assuming that
the algebraic multiplicity of $+1$ is $k$ either,
 Malkin proved \cite{mal} that if the bifurcation
function
$$
  M(h)=\int\limits_0^T \left(\begin{array}{c}
    \left<u_1(\tau,h),g(\tau,x(\tau,\xi(h)),0)\right>\\ ... \\
    \left<u_k(\tau,h),g(\tau,x(\tau,\xi(h)),0)\right>\end{array}\right)
    d\tau
$$
has a simple zero $h_0\in\mathbb{R}^k$, then for any $\varepsilon>0$
sufficiently small system (\ref{mig}) has a unique $T$-periodic
solution $x_\varepsilon$ such that $
   x_\varepsilon(0)\to \xi(h_0)$ as $\varepsilon\to
   0.
$ Here simple zero means that $M(h_0)=0$ and the Jacobian
determinant of $M$ at $h_0$ is nonzero. As usual $\left< \cdot
,\cdot \right>$ denotes the inner product in $\mathbb{R}^n$.
Moreover, Malkin related the asymptotic stability of the solution
$x_\eps$ with the eigenvalues of the Jacobian matrix $DM(h_0).$ The
same result has been proved independently by Loud \cite{loud}.

Malkin's result has been developed by Kopnin \cite{kop} and Loud
\cite{loud} who studied the case when the zero $h_0$ of $M$ is not
necessary  simple  and both  authors obtained conditions for the
existence,  uniqueness and  asymptotic stability of the $T$-periodic
solution $x_\varepsilon$ of  (\ref{mig}) satisfying $
   x_\varepsilon(0)\to \xi(h_0)$ as $\varepsilon\to
   0.$ Other improvements suppress the smoothness
   requirement for $g$ and are due to
Fe\v{c}kan \cite{fec} and Kamenskii-Makarenkov--Nistri \cite{nach},
but only the existence and the calculation of the topological index
of  solutions $x_\eps$ of (\ref{mig}) has been considered.

The analysis of the situation when the algebraic multiplicity of
$+1$ is $n$ goes back to Melnikov \cite{melnikov} and stability
for simple and singular zeros $h_0$ of $M$ has been achieved by
Yagasaki \cite{yagasaki}.

In this paper (Theorems~2 and~3 below are our main results)  we
only assume that (\ref{nnn}) has $k$ linearly independent
$T$-periodic solutions (i.e. the geometric multiplicity of the
multiplier $+1$ of (\ref{nnn}) is $k$. Also, we neither assume
that the zero $h_0$ of $M$ is simple, nor that $g$ is
differentiable. More precisely, we assume that in a small
neighborhood of $h_0$ the topological degree of $M$  is nonzero
and that $M$ is a so-called dilating mapping, while for $g$ we
assume that it is Lipschitz and ``piecewise differentiable" in a
suitable sense. Both assumptions are explicit and weaker than the
Malkin's ones. Note that one of the conditions for $g$, denoted
below by (A9), has its roots in \cite{bd,blm,blmS}. On the other
hand we do not obtain asymptotic stability of the
$T$-periodic solution $x_\varepsilon$ of (\ref{mig}) 
    but we prove its existence and uniqueness, in particular we prove that it is isolated.   In order to study the asymptotic stability
   one can eventually use a result of Kolesov's \cite{kol}.  This result is a
    generalization of the Lyapunov linearization
   stability criterium for Lipschitz systems and it requires the isolateness
   of the $T$-periodic solution $x_\eps.$

In order to prove our main result we need to generalize the
Lyapunov--Schmidt reduction method (see \cite{chic}) for the case of
nonsmooth finite dimensional functions. The application of the
generalized Lyapunov--Schmidt reduction method for proving Theorem~2
is done by following the ideas  of 
\cite{adriana}.  

The paper is organized as follows. In the next section we summarize
our notations. In Section~\ref{secLS} we  generalize the
Lyapunov--Schmidt reduction method for  nonsmooth finite dimensional
functions. In Section~\ref{sec3} we prove Theorem~2 and the main
result of the paper, Theorem~\ref{thm3}. An application of this
theorem is done in
 Section~\ref{sec4}.

\section{Notations}
The following notations will be used throughout this paper.

 Let $n,m,k\in\mathbb{N},$ $k\leq n$, $i\in\mathbb{N}\cup\{0\}$.

 We denote the projection onto the first $k$ coordinates by
$\pi:\mathbb{R}^n \to \mathbb{R}^k$, and the one onto the last $n-k$
coordinates by  $\pi^\perp:\mathbb{R}^n \to \mathbb{R}^{n-k}$.

We denote by $I_{n\times n}$  the identity $n\times n$ matrix, while
$0_{n\times m}$  denotes the null $n\times m$ matrix.

For an $n\times n$ matrix $A$ we denote by $A^*$ the adjoint of $A$,
that in the case  the matrix is real reduces to the transpose.

 We consider a norm in $\mathbb{R}^n$
denoted by $\|\cdot\|$. Let $\Psi$ be an $n\times n$ real matrix.
Then $\|\Psi\|$  denotes the operator norm, i.e. $\|\Psi\|=\sup
_{\|\xi\|=1}\|\Psi \xi\|$.

Let $\xi\in\mathbb{R}^n$ and $\mathcal{Z}\subset\mathbb{R}^n$
compact, then  we denote by
$\rho(\xi,\mathcal{Z})=\min_{\zeta\in\mathcal{Z}}\|\xi-\zeta\|$ the
distance between $\xi$ and $\mathcal{Z}$.

For  $\delta>0$ and $z\in \mathbb{R}^n$ the ball in $\mathbb{R}^n$
centered in $z$  of radius $\delta$ will be denoted by
$B_\delta(z)$.

For a subset $\mathcal{U}\subset\mathbb{R}^n$ we denote by
$\mbox{int}(\mathcal{U})$, $\overline{\mathcal{U}}$ and
$\overline{\mbox{co}} \, \mathcal{U}$ its interior, closure and
closure of the convex hull, respectively.

We denote by $C^i(\mathbb{R}^n,\mathbb{R}^m)$ the set of all
continuous and $i$ times continuously differentiable functions from
$\mathbb{R}^n$ into $\mathbb{R}^m.$

Let $\mathcal{F}\in C^0(\mathbb{R}^n,\mathbb{R}^n)$ be a function
that does not have zeros on the boundary of some open bounded set
$\mathcal{U}\subset\mathbb{R}^n$. Then $d(\mathcal{F},\mathcal{U})$
denotes the Brouwer topological degree of $\mathcal{F}$ on
$\mathcal{U}$ (see \cite{brouwer} or \cite[Ch.~1, \S~3]{krazab}).

For $\mathcal{F}\in C^1(\mathbb{R}^n,\mathbb{R}^m)$,
$D\mathcal{F}$ denotes the Jacobian matrix of $\mathcal{F}$. If
 $\mathbb{R}^n=\mathbb{R}^k\times
\mathbb{R}^{n-k}$ and $\alpha \in \mathbb{R}^k,\beta\in
\mathbb{R}^{n-k}$, then $D_\alpha\mathcal{F}(\cdot,\beta)$ denotes
the Jacobian matrix of $\mathcal{F}(\cdot,\beta).$

For $\mathcal{F}\in C^2(\mathbb{R}^n,\mathbb{R})$, $H\mathcal{F}$
denotes the Hessian matrix of $\mathcal{F}$, i.e. the Jacobian
matrix of the gradient of $\mathcal{F}$.

Let $\delta>0$ be sufficiently small. With $o(\delta)$ we denote a
function of variable $\delta$ such that $o(\delta)/\delta\to 0$ as
$\delta \to 0$, while $O(\delta)$ denotes a function of $\delta$
such that $O(\delta)/\delta$ is bounded as $\delta \to 0$. Here the
functions $o$ and $O$ can depend also on other variables, but the
above properties hold uniformly when these variables lie in a fixed
bounded region.


We say that the function $Q:\mathbb{R}^n\times \mathbb{R}^m\to
\mathbb{R}^m$ is {\it locally uniformly Lipschitz with respect to
its first variable} if for each compact $K\subset \mathbb{R}^n\times
\mathbb{R}^m$ there exists $L>0$ such that $\|
Q(z_1,\lambda)-Q(z_2,\lambda)\|\leq L \| z_1-z_2\|$ for all
$(z_1,\lambda), \, (z_2,\lambda)\in K$.

For any Lebesgue measurable  set $M\subset [0,T]$   we denote by
mes($M$) the Lebesgue measure of $M$.

\section{Lyapunov--Schmidt reduction method for nonsmooth finite dimensional
functions}\label{secLS}

If the continuously differentiable function
$P:\mathbb{R}^n\to\mathbb{R}^n$  vanishes on some set
$\mathcal{Z}\subset\mathbb{R}^n$, then sufficient conditions for the
existence of zeros near $\mathcal{Z}$ of the perturbed function
\begin{equation}\label{a}
F(z,\varepsilon)=P(z)+\varepsilon Q(z,\varepsilon), \quad z\in
\mathbb{R}^n, \,\,\, \eps>0  \mbox{ small enough}
\end{equation}
 can be expressed in terms of  the restrictions to $\mathcal{Z}$ of the functions  $z\mapsto DP(z)$ and
$z\mapsto Q(z,0)$. 
Roughly speaking, this is what is known in the literature as the
Lyapunov--Schmidt reduction method, as it is presented for instance
in \cite{chic,adriana} or \cite[\S 24.8]{krazab}. In these
references it is assumed that $Q$ is a continuously differentiable
function. We show in this section that this last assumption can be
weaken. The following theorem is the main result of this section and
generalizes a theorem of \cite{adriana}.

\begin{thm}\label{thm1} Let $P\in
C^1(\mathbb{R}^n,\mathbb{R}^n)$, let $Q\in
C^0(\mathbb{R}^n\times[0,1],\mathbb{R}^n)$  be locally uniformly
Lipschitz with respect to its first variable, and let
$F:\mathbb{R}^n\times[0,1]\to \mathbb{R}^n$ be given by \eqref{a}.
Assume that $P$ satisfies the following hypotheses.
\begin{itemize}
\item[(A1)] \label{hyp:i}
 There exist an invertible $n\times n$ matrix $S,$ an open ball
 $V\subset\mathbb{R}^k$ with $k\le n$, and a
function
 $\beta_0\in C^1(\overline{V},\mathbb{R}^{n-k})$ such that $P$
vanishes on the set
$\mathcal{Z}=\bigcup\limits_{\alpha\in\overline{V}}
\left\{S\left(\begin{array}{c}\alpha \\
\beta_0(\alpha)\end{array}\right)\right\}$. \item[(A2)]
\label{hyp:ii} For any $z\in \mathcal{Z}$ the  matrix $DP(z)S$ has
in its upper right corner the null $k\times(n-k)$ matrix and in the
lower right corner the $(n-k)\times(n-k)$ matrix $\Delta(z)$ with
${\rm det}\left(\Delta(z)\right)\not=0.$
\end{itemize}
For any $\alpha\in \overline V$ we define
\begin{equation} \label{bif_fun}
  \widehat{Q}(\alpha)= \pi Q\left(S\left(\begin{array}{c}\alpha \\
\beta_0(\alpha)\end{array}\right),0\right).
\end{equation}
Then the following statements hold.
\begin{itemize}
\item[(C1)] For any sequences $(z_m)_{m\geq 1}$ from
$\mathbb{R}^n$ and $(\eps_m)_{m\geq 1}$ from $[0,1]$ such that
$z_m\to z_0\in\mathcal{Z}$,  $\varepsilon_m\to 0$ as $m\to\infty$
and $F(z_{_m},\varepsilon_m)=0$ for any $m\geq 1$, we have
$\widehat{Q}\left(\pi S^{-1}z_0\right)=0$. \item[(C2)] If
$\widehat{Q}:\overline V \to \mathbb{R}^k$ is such that
 $ \widehat{Q}(\alpha)\not=0\,\,{\rm for\ all\ }\alpha\in \partial V$
and
 $ d\left(\widehat{Q},V\right)\not=0$,
then there exists  $\varepsilon_1>0$ sufficiently small such that
for each $\varepsilon\in (0,\eps_1]$ there exists at least one
$z_\varepsilon\in \mathbb{R}^n$ with
$F(z_\varepsilon,\varepsilon)=0$ and
$\rho(z_\varepsilon,\mathcal{Z})\to 0$ as $\varepsilon\to 0.$
\end{itemize}
 In addition we assume that there exists $\alpha_0\in V$ such that
$\,\widehat{Q}(\alpha_0)=0$, $\,\widehat{Q}(\alpha)\neq 0$ for all
$\alpha\in \overline V\setminus \{\alpha_0\}$ and
$d\left(\widehat{Q},V\right)\not=0$, and we denote
$z_0=S\left(\begin{array}{c}\alpha_0
\\\beta_0(\alpha_0)\end{array}\right) $. Moreover we also assume:
\begin{itemize}
\item[(A3)]  $P$ is twice differentiable in the points of
 $\mathcal{Z}$,  and  for each $i\in \overline{1,k}$ and  $z\in {\mathcal{Z}}$
 the Hessian matrix  $HP_i(z)$  is symmetric.
 \item[(A4)] There exists  $\delta_1>0$ and $L_{\widehat Q}>0$ such that
 \[||\widehat Q (\alpha_1)-\widehat Q(\alpha _2)|| \geq
L_{\widehat Q} ||\alpha_1-\alpha_2|| \quad \mbox{ for all }
\alpha_1,\alpha_2\in B_{\delta_1}(\alpha_0). \]
\item[(A5)] For  $\delta>0$ sufficiently small we have that
\begin{eqnarray*}
 \label{sd2}
&&\left\|\pi{Q}\left(z_1+
\zeta,\eps\right)-\pi{Q}\left(z_1,0\right)- \left.\pi{Q}\left(z_2+
\zeta,\eps\right)+\pi{Q}\left(z_2,0\right)\right.\right\| \le\\
\nonumber &&  \frac{o(\delta)}{\delta} \|z_1-z_2\|,
\end{eqnarray*}
 for all  $z_1,z_2\in B_{\delta}(z_0)\cap
\mathcal{Z}$, $\eps\in [0,\delta]$  and $\zeta\in B_{\delta}(0)$.
\end{itemize}
Then the following conclusion holds.
\begin{itemize}
\item[(C3)] There exists $\delta_2>0$ such that for each $\varepsilon \in (0,\eps_1]$  there is exactly one $z_\varepsilon \in B_{\delta_2}(z_0)$
with $F(z_\varepsilon,\varepsilon)=0$. Moreover $z_\varepsilon\to
z_0$ as $\varepsilon\to 0.$
\end{itemize}
\end{thm}

We note that a map  that satisfy (A4) is usually called {\it
dilating map} (cf. \cite{altman}).

 For proving Theorem~\ref{thm1} we shall use the
following version of the
Implicit Function Theorem. 

\begin{lem} \label{implicit}
 Let $P\in
C^1(\mathbb{R}^n,\mathbb{R}^n)$  and let $Q\in
C^0(\mathbb{R}^n\times[0,1],\mathbb{R}^n)$  be locally uniformly
Lipschitz with respect to its first variable. Assume that $P$
satisfies the hypotheses {\em (A1)} and {\em (A2)} of Theorem
\ref{thm1}. Then there exist $\delta_0>0,$ $\eps_0>0$ and a function
$\beta:\overline{V}\times[0,\eps_0]\to\mathbb{R}^{n-k}$ such that
\begin{itemize}
\item[(C4)] $\pi^\bot F\left(S\left(\begin{array}{c}\alpha \\
\beta(\alpha,\eps)\end{array}\right),\eps\right)=0$ for all
$\alpha\in \overline{V}$ and $\eps\in [0,\eps_0]$.
\item[(C5)] $\beta(\alpha,\eps)=\beta_0(\alpha)+\eps\mu(\alpha,\eps)$ where
$\mu:\overline{V}\times(0,\eps_0]\to\mathbb{R}^{n-k}$ is bounded.
Moreover for any $\alpha\in \overline V$ and  $\eps\in[0,\eps_0]$,
$\beta(\alpha,\eps)$ is the only zero of\\
 $\pi^\bot F\left(S\left(\begin{array}{c}\alpha \\
\cdot\end{array}\right),\eps\right)$ in
$B_{\delta_0}(\beta_0(\alpha))$ and $\beta$ is continuous in
$\overline V\times [0,\eps_0]$.
\end{itemize}
In addition if $P$ is twice differentiable in the points of
$\mathcal{Z}$, then
\begin{itemize}
\item[(C6)]
there exists $L_\mu>0$ such that $ \|\mu(\alpha_1,\eps)-\mu(\alpha_2,\eps)\|\le
  L_\mu\|\alpha_1-\alpha_2\|$ for all $\alpha_1,\alpha_2\in\overline{V}$ and $\eps\in (0,\eps_0]$.
\end{itemize}
\end{lem}

\noindent{\bf Proof.} (C4) Let
$\widetilde{F}:\mathbb{R}^k\times\mathbb{R}^{n-k}\times[0,1]\to\mathbb{R}^n$
be defined by
$$
  \widetilde{F}(\alpha,\beta,\varepsilon)=F\left(S\left(\begin{array}{c}
  \alpha\\
  \beta
  \end{array}\right),\varepsilon\right),
$$
and let $\widetilde{P}$, $\widetilde{Q}$ and $\widetilde{\Delta}$ be
defined in a similar way. Now the assumptions (A1) and (A2) become
$\widetilde{P}(\alpha,\beta_0(\alpha))=0$ and, respectively, the
matrix $D\widetilde{P}(\alpha,\beta_0(\alpha))$ has in its upper
right corner the null $k\times(n-k)$ matrix and in the lower right
corner the $(n-k)\times(n-k)$ invertible matrix
$\widetilde{\Delta}{(\alpha,\beta_0(\alpha))}$ for any $\alpha\in
\overline V$. Then
$$
\widetilde{F}(\alpha,\beta_0(\alpha),0)
=0\quad{\rm for\ any\ }\alpha\in\overline{V},
$$
and
\begin{equation}
{\rm
det}\left(D_{\beta}\left(\pi^\bot{\widetilde{F}}\right)(\alpha,\beta_0(\alpha),0)\right)
={\rm
  det}\left(\widetilde{\Delta}{(\alpha,\beta_0(\alpha))}\right)\not=0\quad {\rm for\ any\ } \alpha\in\overline{V}.\label{g1}
\end{equation}
It follows from (\ref{g1}) that there exists a radius $\delta>0$
such that
\begin{equation}\label{g2}
\pi^\bot\widetilde{F}(\alpha,\beta,0)\not=0\quad{\rm for\ any\ }
\beta\in \overline{B_{\delta}({\beta}_0(\alpha))}\backslash \left\{
{\beta}_0(\alpha)\right\},\ \alpha\in \overline{V}.
\end{equation}
 The relations (\ref{g1}) and (\ref{g2}) give (see
\cite[Theorem~6.3]{krazab})
$$
  d(\pi^\bot\widetilde{F}(\alpha,\cdot,0),B_{\delta}({\beta}_0(\alpha)))={\rm sign}\left({\rm
  det}(\widetilde{\Delta}{(\alpha,\beta_0(\alpha)})\right)\not=0,\quad\alpha\in\overline{V}.
$$
Hence, by the continuity of the topological degree with respect to
parameters (using the compactness of $\overline{V}$) there exists
$\eps(\delta)>0$ such that
\begin{equation*}
  d(\pi^\bot\widetilde{F}(\alpha,\cdot,\eps),B_\delta({\beta}_0(\alpha)))\not=0\quad{\rm
  for\ any\ }\eps\in[0,\eps(\delta)],\ \alpha\in\overline{V}.
\end{equation*}
This assures the existence of $\beta(\alpha,\varepsilon)\in
B_\delta({\beta}_0(\alpha))$ such that conclusion (C4) holds with
$\delta_0=\delta$ and $\eps_0=\eps(\delta_0)$.

Without loss of generality we can consider in the sequel that
$\eps(\delta)\to 0$ as $\delta\to 0.$ The value of the radius
$\delta$  eventually may decrease in a finite number of steps during
this proof (consequently, also the value of $\eps(\delta)$).
Sometimes we decrease only the value of $\eps(\delta)$, letting
$\delta$ maintaining its value. Without explicitly mentioning it,
finally, in the statement of the lemma, we replace $\delta_0$ by the
least value of the radius $\delta$ and $\eps_0$ by $\eps(\delta)$.

(C5) Since $P$ and $\beta_0$ are $C^1$ and $\overline{V}$ is
bounded, there exists $\eta
>0$ such that the invertible matrix $\Delta$ defined by (A2) satisfies $||\widetilde{\Delta}({\alpha,\beta_0(\alpha)})||\geq 2\eta$
for all $\alpha\in \overline V$. Using again that $P$ is $C^1$ and
$\widetilde{\Delta}(\alpha,\beta_0(\alpha))=D_{\beta}\left( \pi
^\bot \widetilde{P}\right)(\alpha,\beta_0(\alpha))$ , we obtain that
the radius $\delta>0$ found before at (C4) can be decreased, if
necessary, in such a way that
$||\widetilde{\Delta}({\alpha,\beta_0(\alpha)})-D_{\beta}\left( \pi
^\bot \widetilde{P}\right)(\alpha,\beta)||\leq \eta\,$ for all
$\beta\in B_\delta\left( \beta_0(\alpha)\right)$  and $\alpha\in
\overline V$.  Then $||D_{\beta}\left( \pi ^\bot
\widetilde{P}\right)(\alpha,\beta)||\geq \eta\,$ for all $\beta\in
B_\delta\left( \beta_0(\alpha)\right)$,  $\alpha\in \overline V$.
Applying the generalized Mean Value Theorem (see \cite[Proposition
2.6.5]{clark}) to the function $\pi ^\bot
\widetilde{P}(\alpha,\cdot)$, we obtain
\begin{equation} \label{Pinv}
||\pi^\bot \widetilde{P}(\alpha,\beta_1)-\pi^\bot \widetilde{P}(\alpha,\beta_2)||\geq \eta ||\beta_1-\beta_2||,\quad
\beta_1,\beta_2\in B_\delta(\beta_0(\alpha)),\,\,\alpha
\in \overline V.
\end{equation}
We take $M_Q>0$ such that
$||\widetilde{Q}(\alpha,\beta(\alpha,\eps),\eps)||\leq M_Q$ for all
$\alpha\in \overline V$ and $\eps\in[0,\eps_0]$. Using \eqref{Pinv}
we obtain for all $\alpha\in \overline V$ and
$\eps\in[0,\eps(\delta)]$
\begin{eqnarray*} 0&=&||\pi^\bot
\widetilde{P}(\alpha,\beta(\alpha,\eps))-\pi^\bot \widetilde{P}(\alpha,\beta_0(\alpha))+\eps
\pi^\bot
\widetilde{Q}(\alpha,\beta(\alpha,\eps),\eps)||\\
&\geq& \eta ||\beta(\alpha,\eps)-\beta_0(\alpha)||-\eps M_Q\,.
\end{eqnarray*}
From these last relations, denoting $m=M_Q/\eta$, we deduce that
\begin{equation} \label{mu_bdd}
\|\mu(\alpha,\eps)\|\leq m \quad \mbox{ for all }\alpha\in \overline
V, \, \eps\in (0,\eps(\delta)].
\end{equation}
We choose $L_Q>0$ such that
\begin{equation} \label{Q_Lip} || \widetilde{Q}(\alpha_2,\beta_2,\eps)-
\widetilde{Q}(\alpha_1,\beta_1,\eps)||\leq
L_Q\left(||\alpha_2-\alpha_1||+ ||\beta_2-\beta_1||\right)\, ,
\end{equation} for all $\beta_1,\beta_2\in B_{\delta_0}(\beta_0(\overline V))$,
$\alpha_1,\alpha_2\in \overline V,\eps\in [0,\eps_0].$ We decrease
$\delta>0$ in such a way that $\eta-\eps L_Q>0$ for any
$\eps\in[0,\eps(\delta)].$

 Let $\alpha\in \overline V$, $\eps\in[0,\eps(\delta)]$ and assume  that $\beta(\alpha,\eps)$ and
$\beta_2$ are two zeros of $\pi^\bot F\left(S\left(\begin{array}{c}\alpha \\
\cdot\end{array}\right),\eps\right)$ in $B_\delta(\beta_0(\alpha))$.
Taking into account (\ref{Pinv}) and \eqref{Q_Lip}, we obtain
\begin{eqnarray*}
0&=&||\pi^\bot \widetilde{P}(\alpha,\beta_2)-\pi^\bot
\widetilde{P}(\alpha,\beta(\alpha,\eps))+\\
&~& \eps \pi^\bot \widetilde{Q}(\alpha,\beta_2,\eps)-\eps \pi^\bot
\widetilde{Q}(\alpha,\beta(\alpha,\eps),\eps)||\\
&\geq& (\eta-\eps L_Q) ||\beta_2-\beta(\alpha,\eps)||.
\end{eqnarray*}
Since $\eta-\eps L_Q>0$ for any $\eps\in[0,\eps(\delta)]$ we deduce
from this last relation that $\beta_2$ and $\beta(\alpha,\eps)$ must
coincide.

\medskip

We prove in the sequel the continuity of the function
$\beta:\overline V \times [0,\eps(\delta)]\to \rkn$. Let
$(\alpha_1,\eps_1)\in \overline V \times [0,\eps(\delta)]$ be fixed
and $(\alpha,\eps)\in \overline V \times [0,\eps(\delta)]$ be in a
small neighborhood of $(\alpha_1,\eps_1)$. Consider $L_P>0$ such
that $||\widetilde{P}(\alpha_1,\beta)-\widetilde P
(\alpha,\beta)||\leq L_P||\alpha_1-\alpha||$ for all
$\alpha_1,\alpha\in \overline V$ and $\beta\in
B_{\delta_0}(\beta_0(\overline V))$.   We diminish $\eps(\delta)>0$,
if necessary, and we consider $\alpha$ so close to $\alpha_1$ that
$\beta(\alpha,\eps)\in B_\delta(\beta_0(\alpha_1))$.
Then using (\ref{Pinv}) and (\ref{Q_Lip}) we obtain
\begin{eqnarray*}
0&=&||\pi^\bot \widetilde{P}(\alpha_1,\beta(\alpha_1,\eps_1))-\pi^\bot
\widetilde{P}(\alpha,\beta(\alpha,\eps))+\\
&~& \eps_1 \pi^\bot
\widetilde{Q}(\alpha_1,\beta(\alpha_1,\eps_1),\eps_1)-\eps \pi^\bot
\widetilde{Q}(\alpha,\beta(\alpha,\eps),\eps)|| \\
&\geq& \eta
||\beta(\alpha_1,\eps_1)-\beta(\alpha,\eps)||-L_P||\alpha_1-\alpha||-\\
&~& ||\eps_1\pi^\bot
\widetilde{Q}(\alpha_1,\beta(\alpha_1,\eps_1),\eps_1)-\eps \pi^\bot
\widetilde{Q}(\alpha,\beta(\alpha,\eps),\eps) ||
\end{eqnarray*}
and
\begin{eqnarray*}
&& -||\eps_1\pi^\bot
\widetilde{Q}(\alpha_1,\beta(\alpha_1,\eps_1),\eps_1)-\eps
\pi^\bot
\widetilde{Q}(\alpha,\beta(\alpha,\eps),\eps) ||  \\
&\ge & -\eps_1 L_Q ||\alpha_1-\alpha||-\eps_1
L_Q||\beta(\alpha_1,\eps_1)-\beta(\alpha,\eps)||-\\
&~&||\eps_1\pi^\bot
\widetilde{Q}(\alpha,\beta(\alpha,\eps),\eps_1)-\eps \pi^\bot
\widetilde{Q}(\alpha,\beta(\alpha,\eps),\eps)||.
\end{eqnarray*}
Combining these last two relations we obtain
\begin{eqnarray*}
&&(\eta -\eps_1
L_Q)||\beta(\alpha_1,\eps_1)-\beta(\alpha,\eps)||\leq (L_P+\eps_1
L_Q) ||\alpha_1-\alpha||+\\
&~&||\eps_1\pi^\bot
\widetilde{Q}(\alpha,\beta(\alpha,\eps),\eps_1)-\eps \pi^\bot
\widetilde{Q}(\alpha,\beta(\alpha,\eps),\eps)||\, ,
\end{eqnarray*}
from where it follows easily that $\beta(\alpha,\eps)\to \beta (\alpha_1,\eps_1)$ when
$(\alpha,\eps)\to (\alpha_1,\eps_1)$. \\

(C6)  We define $\Phi(\alpha,\xi)=\pi^\bot
\widetilde{P}(\alpha,\beta_0(\alpha)+\xi)$ for all
$\alpha\in\overline V$ and $\xi\in\mathbb{R}^{n-k}$. From
(\ref{Pinv}) we have that
\begin{equation} \label{Phiinv}
||\Phi(\alpha,\xi_1)-\Phi(\alpha,\xi_2)||\geq \eta ||\xi_1-\xi_2||
\mbox{ for all } \alpha\in \overline V, \,\, \xi_1,\xi_2\in
B_\delta(0).
\end{equation}
 Since $\widetilde P(\alpha,\beta_0(\alpha))=0$ for all $\alpha\in
\overline V$, we have that $\Phi(\alpha,\xi)=\pi^\bot
\widetilde{P}(\alpha,\beta_0(\alpha)+\xi)-\pi^\bot
\widetilde{P}(\alpha,\beta_0(\alpha))$ and that
\begin{eqnarray*}
D_\alpha  \Phi(\alpha,\xi)&=&D_\alpha \left(\pi^\bot
\widetilde{P}\right)(\alpha,\beta_0(\alpha)+\xi) - D_\alpha \left(\pi^\bot
\widetilde{P}\right)(\alpha,\beta_0(\alpha)) +\\
&&\left[ D_\beta \left( \pi^\bot \widetilde{P}\right)(\alpha,\beta_0(\alpha)+\xi)-D_\beta \left(
\pi^\bot \widetilde{P}\right)(\alpha,\beta_0(\alpha))\right]D\beta_0(\alpha).
\end{eqnarray*}
From this expression, using that $ \widetilde{P}$ is twice
differentiable in $(\alpha,\beta_0(\alpha))$ and $\beta_0$ is $C^1$,
we obtain  for some $L_\Phi>0$ that the radius $\delta$ can be
eventually decreased in a such way that
$$||D_\alpha \Phi(\alpha,\xi)||\leq
L_{\Phi}||\xi|| \,\mbox{ for all }\alpha\in \overline V,\, \xi\in
B_\delta(0).$$
 Hence using  the Mean Value Inequality we have
\begin{equation} \label{Philip}
||\Phi(\alpha_1,\xi)-\Phi(\alpha_2,\xi)||\leq L_{\phi}||\xi||\cdot
||\alpha_1-\alpha_2|| \, \mbox{ for all }\alpha_1,\alpha_2\in
\overline V,\, \xi\in B_\delta(0).
\end{equation}
Now we use (\ref{Phiinv}) with $\xi_1=\eps \mu(\alpha_1,\eps)$,
$\xi_2=\eps\mu(\alpha_2,\eps)$  diminishing $\eps(\delta)$, if
necessary, in order that $\xi_1,\xi_2\in B_\delta(0)$ for all
$\alpha_1,\alpha_2\in\overline V$ and $\eps\in(0,\eps(\delta)]$.
Using also (C5), (\ref{mu_bdd}) and
 (\ref{Philip})  we obtain
\begin{equation} \label{eq:despreP} \begin{array}{ll}
&||\pi^\bot \widetilde{P}(\alpha_1,\beta(\alpha_1,\eps))-\pi^\bot
\widetilde{P}(\alpha_2,\beta(\alpha_2,\eps))||=
 ||\Phi(\alpha_1,\xi_1)-\Phi(\alpha_2,\xi_2)|| \\
&\geq \eta ||\xi_1-\xi_2||- L_\Phi ||\xi_1||\cdot
||\alpha_1-\alpha_2||
\\ &\geq \eta\eps||\mu(\alpha_1,\eps)-\mu(\alpha_2,\eps)||-L_\Phi m \eps
||\alpha_1-\alpha_2||\ ,
\end{array}
\end{equation}
for all $\alpha_1,\alpha_2\in \overline V$ and $\eps\in
(0,\eps(\delta)]$. Also using (\ref{Q_Lip}) we have
\begin{equation} \label{eq:despreQ}
\begin{array}{ll}
||\pi^\bot \widetilde{Q}(\alpha_1,\beta(\alpha_1,\eps),\eps)- \pi^\bot
\widetilde{Q}(\alpha_2,\beta(\alpha_2,\eps),\eps) || \leq \\
\leq  \eps L_Q||\mu(\alpha_1,\eps)-\mu(\alpha_2,\eps)||+
L_Q(1+L_{\beta_0}) ||\alpha_1-\alpha_2||\ ,
\end{array}
\end{equation}
for all $\alpha_1,\alpha_2\in \overline V$ and $\eps\in
(0,\eps(\delta)]$, where  $L_{\beta_0}$ is the Lipschitz constant of
$\beta_0$ in $\overline V$. By definition of $\beta(\alpha,\eps)$ we
have $\,\pi^\bot
\widetilde{P}(\alpha_i,\beta(\alpha_i,\eps))+\eps\pi^\bot
\widetilde{Q}(\alpha_i,\beta(\alpha_i,\eps),\eps)=0$ for $i\in
\overline{1,2}$. Using (\ref{eq:despreP}) and (\ref{eq:despreQ}) we
obtain
\[ 0\geq \eps [\eta  -\eps L_Q ]\cdot ||\mu(\alpha_1,\eps)-\mu(\alpha_2,\eps)||-\eps [L_\Phi m
+L_Q(1+L_{\beta_0})]\cdot ||\alpha_1-\alpha_2||\ ,
\] 
for all $\alpha_1,\alpha_2\in \overline V$ and $\eps\in
(0,\eps(\delta)]$. Therefore $\mu :\overline
V\times(0,\eps(\delta)]\to \rkn$ satisfies (C6) with $L_\mu =[L_\Phi
m +L_Q(1+L_{\beta_0})]/[\eta -\eps(\delta) L_Q ]$. Hence all the
conclusions hold with $\delta_0=\delta$ and $\eps_0=\eps(\delta)$.

\qed

We remark that (C4) and the uniqueness part of (C5) can be obtained
by means of the Lipschitz generalization of the Inverse Function
Theorem (see e.g. \cite[Theorem 5.3.8]{krantz}), but we provide a
different proof because the inequalities (\ref{Pinv}) and
(\ref{mu_bdd}) are used for proving the rest of (C5) and (C6).

\

\noindent{\bf Proof of Theorem~\ref{thm1}.} Let $\delta_0$,
$\eps_0$, $\beta(\alpha,\eps)$ and $\mu(\alpha,\eps)$ be as in
 Lemma~\ref{implicit}. We consider the notations
$\widetilde F,\widetilde P$ and $\widetilde Q$ like in  the proof of
Lemma \ref{implicit}.

\medskip

(C1) Let the sequences $(z_m)_{m\geq 1}$ from $\mathbb{R}^n$ and
$(\eps_m)_{m\geq 1}$ from $[0,1]$ be such that $z_m\to z_0\in
{\mathcal{Z}}$, $\eps_m\to 0$ as $m\to\infty$ and $F(z_m,\eps_m)=0$
for any $m\geq 1$. We define $\alpha_0\in \mathbb{R}^k$, the
sequences $(\alpha_m)_{m\geq 1}$ from $\mathbb{R}^k$ and
$(\beta_m)_{m\geq 1}$ from $\mathbb{R}^{n-k}$ by  $z_0=S\left(\begin{array}{c}\alpha_0\\
\beta_0(\alpha_0)\end{array}\right)$ and $z_m=S\left(\begin{array}{c}\alpha_m\\
\beta_m\end{array}\right)$.  Then we have that $\alpha_0=\lim
\limits_{m\to\infty}\alpha_m$,
$\beta_0(\alpha_0)=\lim \limits_{m\to\infty}\beta_m$  
and there exists
$m_0\in\mathbb{N}$ such that $\beta_m\in
B_{\delta_0}(\beta_0(\alpha_m))$ and $\eps_m\in [0,\eps_0]$ for all $m\geq m_0$. 
Therefore, since $F(z_m,\eps_m)=0$, Lemma~\ref{implicit} implies
$\beta_m=\beta(\alpha_m,\eps_m)$ for any $m\geq m_0.$ Since
$\pi\widetilde{P}(\alpha_m,\beta_0(\alpha_m))=0$ and $D_\beta
(\pi\widetilde{P})(\alpha_m,\beta_0(\alpha_m))=0$, we obtain that
$\lim \limits_{m\to
\infty}\dfrac{1}{\eps_m}\pi\widetilde{P}(\alpha_m,\beta(\alpha_m,\eps_m))=0.$
   Hence
\begin{eqnarray*}
  0&=&\lim \limits_{m\to \infty}\frac{1}{\eps_m} \pi
  \widetilde{F}(\alpha_m,\beta(\alpha_m,\eps_m),\eps_m)\\
  &=&\lim \limits_{m\to \infty}\left[\frac{1}{\eps_m}
  \pi\widetilde{P}(\alpha_m,\beta(\alpha_m,\eps_m))+\pi\widetilde{Q}(\alpha_m,\beta(\alpha_m,\eps_m),\eps_m)\right]=
\widehat{Q}(\alpha_0) \end{eqnarray*} from where (C1)  follows.

\medskip

(C2) Using  (C4) of Lemma~\ref{implicit}, we note that it is enough
to prove the existence  of at least one zero in $V$ of the function
$\alpha\mapsto  \pi \widetilde{F}(\alpha,\beta(\alpha,\eps),\eps)$
for each $\eps\in(0,\eps_1]$ where $\eps_1$ with $0<\eps_1\leq
\eps_0$ has to be found. This will follow from the claim that the
Brouwer topological degree $d\left(\frac{1}{\eps} \pi
  \widetilde{F}(\cdot,\beta(\cdot,\eps),\eps),V\right)\neq 0$ for $\eps\in
  (0,\eps_1]$.
 Now we prove this claim. Since
$\beta(\alpha,\eps)=\beta_0(\alpha)+\eps \mu(\alpha,\eps)$ with
$\mu:\overline V\times (0,\eps_0]\to \mathbb{R}^{n-k}$ a bounded
function, $\pi\widetilde{P}(\alpha,\beta_0(\alpha))=0$ and $D_\beta
(\pi\widetilde{P})(\alpha,\beta_0(\alpha))=0$, we have
\[ \lim_{\eps\to 0}\frac{1}{\eps}\pi\widetilde{P}(\alpha,\beta(\alpha,\eps))=0.\]
   Therefore
$$
  \lim_{\eps\to 0}\frac{1}{\eps} \pi
  \widetilde{F}(\alpha,\beta(\alpha,\eps),\eps)=\lim_{\eps\to 0}\left[\frac{1}{\eps}
  \pi\widetilde{P}(\alpha,\beta(\alpha,\eps))+\pi\widetilde{Q}(\alpha,\beta(\alpha,\eps),\eps)\right]=
\widehat{Q}(\alpha).
$$
Using the continuity of the Brouwer degree with respect to the
parameter $\eps$, and taking into account that, by hypothesis,
$d(\widehat{Q},V)\not=0$,  for each $\eps\in (0,\eps_1]$ there
exists $\eps_1>0$ sufficiently small such that
$$
  d\left(\frac{1}{\eps}\pi\widetilde{F}(\cdot,\beta(\cdot,\eps),\eps),V\right)=d(\widehat{Q},V)\not=0.
$$
Hence the claim is proved. Then for each $\eps \in (0,\eps_1]$
 there exists $\alpha_\eps \in V$ such that
$\pi\widetilde{F}(\alpha_\eps,\beta(\alpha_\eps,\eps),\eps)=0$ and,
moreover, using also (C4) of Lemma~\ref{implicit}, we have that
$\widetilde{F}(\alpha_\eps,\beta(\alpha_\eps,\eps),\eps)=0$. Denoting $z_\eps=S\left(\begin{array}{c} \alpha_\eps \\
\beta(\alpha_\eps,\eps) \end{array} \right)$ we have  that
$F(z_\eps,\eps)=0$. From the definitions of $z_\eps$ and
$\mathcal{Z}$, and the continuity of $\beta$, it follows easily that
$\rho(z_\eps,\mathcal{Z})\to 0$ as $\eps\to 0$.

\medskip

(C3)  Since $\alpha_0\in V$ is an isolated zero of $\widehat Q$,
applying the topological degree arguments like in (C2) for $V$ that
shrinks
  to $\{\alpha_0\}$, we obtain the existence of $\alpha_{\eps}$ such that $\alpha_{\eps}\to \alpha_0$ as $\eps\to
  0$,
  and $\pi\widetilde{F}(\alpha_\eps,\beta(\alpha_\eps,\eps),\eps)=0$ for any $\eps\in
  (0,\eps_1]$.
  Hence $z_\eps=S\left(\begin{array}{c} \alpha_\eps \\
\beta(\alpha_\eps,\eps) \end{array} \right)$ and
$z_0=S\left(\begin{array}{c} \alpha_0 \\ \beta_0(\alpha_0)
\end{array} \right)\in \mathcal{Z}$ are  such that
$F(z_\eps,\eps)=0$ and $z_\eps \to z_0$ as $\eps\to 0$.

In order to prove that $z_\eps$ is the unique zero of
$F(\cdot,\eps)$ in a neighborhood of $z_0$, we define
$$r_1(\alpha,\eps)=\frac{1}
  {\eps}\pi\widetilde{P}(\alpha,\beta(\alpha,\eps)), \quad
  r_2(\alpha,\eps)=\pi\widetilde{Q}(\alpha,\beta(\alpha,\eps),\eps)-\pi\widetilde{Q}(\alpha,\beta_0(\alpha),0),$$
for all $\alpha\in\overline V$ and $\eps\in(0,\eps_1]$, and we study
the Lipschitz properties with respect to $\alpha$ of these two
functions.

Since $\widetilde{P}(\alpha,\beta_0(\alpha))=0$ for all $\alpha\in
\overline V$, by taking the derivative with respect to $\alpha$ we
obtain
\begin{equation}\label{tmp}
  D_\alpha \left(\pi \widetilde{P}\right)(\alpha,\beta_0(\alpha))+D_\beta \left(\pi \widetilde{P}\right)
  (\alpha,\beta_0(\alpha))D \beta_0(\alpha)=0\quad{\rm
  for\ all\ }\alpha\in \overline V.
\end{equation}

Assumption (A2) assures that
$D_\beta \left(\pi \widetilde{P}\right)(\alpha,\beta_0(\alpha))
=0$
 for all $\alpha\in \overline V$. Taking the derivative with respect to $\alpha$, we have
\begin{equation} \label{secondorder}
D_{\beta\alpha}\left(\pi
{\widetilde{P}}\right)(\alpha,\beta_0(\alpha))+D_{\beta\beta}\left(\pi
{\widetilde{P}}\right)(\alpha,\beta_0(\alpha))D\beta_0(\alpha)=0\quad{\rm
for\ any\ }\alpha\in \overline V.
\end{equation}
For any $\alpha\in\overline V$ and $\xi\in \mathbb{R}^{n-k}$ we
define $\Phi(\alpha,\xi)=\pi
\widetilde{P}(\alpha,\beta_0(\alpha)+\xi)$. Taking into account the
relations (\ref{tmp}) and
(\ref{secondorder}) and that, by hypothesis (A3) we have that\\
$D_{\beta\alpha}\left(\pi
{\widetilde{P}}\right)(\alpha,\beta_0(\alpha))=D_{\alpha\beta}\left(\pi
{\widetilde{P}}\right)(\alpha,\beta_0(\alpha))$, we obtain
\begin{eqnarray*}
D_{\alpha}\Phi (\alpha,\xi)&=& D_\alpha \left(\pi
\widetilde{P}\right)(\alpha,\beta_0(\alpha)+\xi)+D_\beta \left(\pi \widetilde{P}\right)
  (\alpha,\beta_0(\alpha)+\xi)D \beta_0(\alpha)-\\
&~& D_\alpha \left(\pi
\widetilde{P}\right)(\alpha,\beta_0(\alpha))-D_\beta \left(\pi
\widetilde{P}\right)
  (\alpha,\beta_0(\alpha))D \beta_0(\alpha)-\\
&~&D_{\alpha\beta}\left(\pi
{\widetilde{P}}\right)(\alpha,\beta_0(\alpha))\xi-D_{\beta\beta}\left(\pi
{\widetilde{P}}\right)(\alpha,\beta_0(\alpha))D\beta_0(\alpha)\xi.
\end{eqnarray*}
From this last equality, using that $D_\alpha \left(\pi
\widetilde{P}\right)$ and, respectively, $D_\beta \left(\pi
\widetilde{P}\right)$ are differentiable at
$(\alpha,\beta_0(\alpha))$, we deduce that $D_{\alpha}\Phi
(\alpha,\xi)=o(\xi)$ for all $\alpha\in\overline V$ and $\xi\in
\mathbb{R}^{n-k}$ with $||\xi||$ sufficiently small. Hence the mean
value inequality assures that
$$||\Phi(\alpha_1,\xi)-\Phi(\alpha_2,\xi)||\leq
o(\xi)||\alpha_1-\alpha_2|| \quad \mbox{ for all }
\alpha_1,\alpha_2\in\overline V.$$  In the last inequality we
replace $\xi = \eps \mu (\alpha_1,\eps)$ (where $\mu$ is given by
Lemma \ref{implicit}).  We use that $D_\xi \Phi(\alpha,0)=D_\beta
\pi\widetilde{P}(\alpha,\beta_0(\alpha))=0$ for any
$\alpha\in\overline V,$ and that $\mu$ is Lipschitz with respect to
$\alpha \in\overline V$. Then we obtain, considering that $\eps_1$
is small enough,  for all $\eps\in (0,\eps_1]$
$$||\Phi(\alpha_1,\eps\mu(\alpha_1,\eps))-\Phi(\alpha_2,\eps \mu(\alpha_2,\eps))||\leq
o(\eps)||\alpha_1-\alpha_2|| \quad \mbox{ for all }
\alpha_1,\alpha_2\in V.$$ Now coming back to our notations and
recalling that $\beta(\alpha,\eps)=\beta_0(\alpha)+\eps \mu
(\alpha,\eps)$, we obtain for $\eps\in (0,\eps_1]$
\begin{equation} \label{lr1} ||r_1(\alpha_1,\eps)-r_1(\alpha_2,\eps)||\leq
\frac{o(\eps)}{\eps}||\alpha_1-\alpha_2|| \quad \mbox{ for all }
\alpha_1,\alpha_2\in\overline V.\end{equation}

We will prove that a similar relation holds for the function $r_2$.
First we note that the hypothesis (A5) and the fact that $Q$ is
locally uniformly Lipschitz with respect to the first variable imply
that 
 \begin{equation}\label{sd3}
\begin{array}{c}
\left\|\pi{Q}\left(z_1+
\zeta_1,\eps\right)-\pi{Q}\left(z_1,0\right)-
\left.\pi{Q}\left(z_2+ \zeta_2,\eps\right)+\pi{Q}\left(z_2,0\right)\right.\right\|\le\\

 \frac{o(\delta)}{\delta} \|z_1-z_2\|+L_Q\|\zeta_1-\zeta_2\|,
 \end{array}
\end{equation}
 for all  $z_1,z_2\in B_{\delta}(z_0)\cap
\mathcal{Z}$, $\eps\in [0,\delta]$  and $\zeta_1,\zeta_2\in
B_{\delta}(0)$. 
 We diminish
$\delta_1>0$ given in (A4) and $\eps_1>0$ in such a way
that $\delta_1\leq \delta$, $\eps_1\leq \delta$, $S\left(\begin{array}{c}\alpha\\
\beta_0(\alpha)\end{array}\right)\in B_{\delta}(z_0)$ and $S\left(\begin{array}{c}0_{k\times 1}\\
\eps\mu(\alpha,\eps)\end{array}\right)\in B_{\delta}(0)$ for any
$\alpha\in B_{\delta_1}(\alpha_0),$ $\eps\in(0,\eps_1].$ Replacing
$z_i=S\left(\begin{array}{c}
  \alpha_i\\
  \beta_0(\alpha_i)
  \end{array}\right)$,
$\zeta_i= S\left(\begin{array}{c}
  0_{k\times1}\\
  \eps \mu(\alpha_i,\eps)
  \end{array}\right)$, $i\in \overline{1,2}$ in \eqref{sd3} we obtain  that
 \begin{eqnarray*}
&&||r_2(\alpha_1,\eps)-r_2(\alpha_2,\eps)||\leq \\
&& \leq \frac{o(\delta)}{\delta} \left(\, ||\alpha_1-\alpha_2|| +
||\beta_0(\alpha_1)-\beta_0(\alpha_2)|| \,\right) + \eps L_Q
||\mu(\alpha_1,\eps)-\mu (\alpha_2,\eps)||\ ,
 \end{eqnarray*}
 for all $\alpha_1,\alpha_2\in B_{\delta_1}(\alpha_0)$ and $\eps\in(0,\eps_1]$.
By hypothesis, $\beta_0$ is $C^1$ in $\overline V$ and, by
Lemma~\ref{implicit} (conclusion (C6)), $(\alpha,\eps)\mapsto
\mu(\alpha,\eps)$ is Lipschitz with respect to $\alpha\in\overline
V$ (with a Lipschitz constant that does not
depend on $\eps$). Hence  for 
 $\delta_1,\eps_1\leq \delta$ small enough, 
\begin{equation} \label{lr2}
||r_2(\alpha_1,\eps)-r_2(\alpha_2,\eps)||\leq
\frac{o(\delta)}{\delta}||\alpha_1-\alpha_2||, \quad
\alpha_1,\alpha_2\in B_{\delta_1}(\alpha_0), \, \eps\in
(0,\eps_1].\end{equation} Therefore we have proved that $r_1$ and
$r_2$ satisfy the Lipschitz conditions \eqref{lr1} and,
respectively, \eqref{lr2}. In what follows we define some constant
$\delta_2>0$, and after we prove that it is the one that satisfies
the requirements of (C3).


We diminish $\delta_1>0$ in such a way that there exists
$\delta_3>0$ such that $\delta_3\leq \delta_0$ and
$B_{\delta_3}(\beta_0(\alpha_0))\subset \bigcap\limits_{\alpha\in
B_{\delta_1}(\alpha_0)}B_{\delta_0}(\beta_0(\alpha))$. We choose
$\delta_2>0$ so small that $S^{-1}(B_{\delta_2}(z_0))\subset
B_{\delta_1}(\alpha_0)\times B_{\delta_3}(\beta_0(\alpha_0))$. We
diminish $\eps_1>0$, if necessary,  such that   $z_\eps \in
B_{\delta_2}(z_0)$ for any $\eps\in(0,\eps_1].$ For any
$\eps\in(0,\eps_1]$ we claim that
 $z_\eps$ is the only zero of $F(\cdot,\eps)$ in
$B_{\delta_2}(z_0).$ Assume by contradiction that there exists
$\eps_2\in(0,\eps_1]$ such that $z_{\eps_2}$ and $z_2$ are two
different zeros of $F(\cdot,\eps_2)$ in $B_{\delta_2}(z_0).$
 Denoting $\alpha_2=\pi S^{-1}z_2$ and  $\beta_2=\pi^{\bot}
 S^{-1}z_2$ we have that $\beta_2\in
 B_{\delta_0}(\beta_0(\alpha_2))$. By (C5) of Lemma~\ref{implicit},
 since $\beta_2$ is a zero of $\pi ^{\bot} F \left( S \left( \begin{array}{c} \alpha_2 \\ \cdot \end{array}\right),\eps_2
 \right)$ (using the notations introduced before, $\pi ^{\bot}\widetilde{F}(\alpha_2,\cdot,\eps_2)$),  we must have that $\beta_2=\beta(\alpha_2,\eps_2)$.
 Therefore $\alpha_{\eps_2}$ and $\alpha_2$ are two different zeros of $\pi
\widetilde{F}(\cdot,\beta(\cdot,\eps_2),\eps_2)$ in
$B_{\delta_1}(\alpha_0).$ We have the identity
$$
\frac{1}{\eps}\pi\widetilde{F}(\alpha,\beta(\alpha,\eps),\eps)=\widehat{Q}(\alpha)+r_1(\alpha,\eps)+r_2(\alpha,\eps)
\mbox{ for all } \alpha\in \overline V ,\,\,\, \eps\in (0,\eps_1].$$

We denote $r(\alpha,\eps)=r_1(\alpha,\eps)+r_2(\alpha,\eps)$.
  Then assumption (A4), properties (\ref{lr1}) and
(\ref{lr2}) give
$$0=||\widehat{Q}(\alpha_{\eps_2})-\widehat{Q}(\alpha_2)+r(\alpha_{\eps_2},\eps_2)-r(\alpha_2,\eps_2)||\geq
(L_{\widehat{Q}}-o(\eps_2)/{\eps_2}-o(\delta)/\delta)||\alpha_{\eps_2}-\alpha_2||.$$
Since $\eps_1>0$ and $\delta>0$ are sufficiently small and
$0<\eps_2\leq \eps_1$, the constant
$(L_{\widehat{Q}}-o(\eps_2)/{\eps_2}-o(\delta)/\delta)$ must be
positive and, consequently,  $\alpha_{\eps_2}$ and $\alpha_2$ must
coincide. Hence also $z_{\eps_2}$ and $z_2$ must coincide and we
conclude the proof.

\qed

\section{A generalization of Malkin's result on the existence of $T$-periodic
solutions for $T$-periodically perturbed differential equations when
the perturbation is nonsmooth}\label{sec3}

In this section we consider the problem of  existence and uniqueness
of $T$-periodic solutions for the $T$-periodic differential system
\begin{equation}
\label{ps} \dot{x}= f(t,x)+\eps g(t,x,\eps)\, ,
\end{equation}
where $f\in C^2(\mathbb{R}\times\mathbb{R}^n, \mathbb{R}^n)$ and
$g\in C^0(\mathbb{R}\times\mathbb{R}^n\times[0,1],\mathbb{R}^n)$ are
$T$-periodic in the first variable and  $g$ is locally uniformly
Lipschitz with respect to its second variable. For $z\in \rn$ we
denote by $x(\cdot,z,\eps)$ the solution of \eqref{ps} such that
$x(0,z,\eps)=z$. We consider the situation when the unperturbed
system
\begin{equation}\label{np}
\dot{x}= f(t,x)\, ,
\end{equation}
has a non--degenerate (in a sense that will be precised below)
family of $T$-periodic solutions.
The main tool for the proof of our main result is Theorem
\ref{thm1}. We will show that the assumptions of Theorem~\ref{thm1}
can be expressed in terms of the function $g$ and of the solutions
of the linear differential system
\begin{equation} \label{ls} \dot{y}=D_xf(t,x(t,z,0))y. \end{equation}
Indeed we have the following theorem, which generalizes a related
result by Roseau and improves it with the uniqueness of the
periodic solution. The above mentioned result by Roseau is proved
in a shorter way in \cite{adriana}. Here we will use the same main
ideas from \cite{adriana} to prove the next result.

\begin{thm}\label{thm2} Assume that $f\in C^2(\mathbb{R}\times\mathbb{R}^n, \mathbb{R}^n)$ and
 $g\in C^0(\mathbb{R}\times\mathbb{R}^n\times[0,1],\mathbb{R}^n)$ are $T$-periodic in the first
variable, and  that $g$ is locally uniformly Lipschitz with respect
to the second variable. Suppose that the unperturbed system
\eqref{np}
 satisfies the following conditions.
\begin{itemize}
\item[(A6)] There exist an invertible $n\times n$ real matrix
$S$, an open ball  
 $V\subset {\mathbb{R}^k}$ with $k\leq n$, and
  a $C^1$
function $\beta_0:\overline V\to {\mathbb{R}^{n-k}}$ such that any
point of the set $\mathcal{Z}=\bigcup\limits_{\alpha\in \overline
V}\left\{S\left(\begin{array}{c} \alpha\\
\beta_0(\alpha)\end{array}\right)\right\}$ is the initial condition
of a $T$-periodic solution of \eqref{np}.

\item[(A7)] For each $z\in \mathcal{Z}$ there exists a
fundamental matrix solution $Y(\cdot,z)$ of (\ref{ls}) such that
***$Y(0,z)$ is $C^1$ with respect to $z$ and *** such that the matrix
$\left(Y^{-1}(0,z)-Y^{-1}(T,z)\right)S$ has in the upper right
corner the null $k\times (n-k)$ matrix, while in the
 lower right corner has the $(n-k)\times (n-k)$ matrix $\Delta(z)$ with $\det(\Delta(z))\neq
 0$.
\end{itemize}
\noindent We define the function $G:\overline V \to {\mathbb{R}^k}$  by
\begin{equation} \nonumber 
G(\alpha)=\pi \int _0^T
Y^{-1}\left(t,S\left(\begin{array}{c} \alpha\\
\beta_0(\alpha)\end{array}\right)\right)g\left(t,x\left(t,S\left(\begin{array}{c} \alpha\\
\beta_0(\alpha)\end{array}\right) ,0 \right),0\right) dt.
\end{equation}
Then the following statements hold.
\begin{itemize}
\item[(C7)] For any  sequences $(\varphi_m)_{m\geq 1}$ from $C^0(\mathbb{R}, \mathbb{R}^n)$ and $(\eps_m)_{m\geq 1}$
from $[0,1]$ such that $\varphi_m(0)\to z_0\in \mathcal{Z}$,
$\eps_m\to 0$ as $m\to \infty$ and $\varphi_m$ is a $T$-periodic
solution of \eqref{ps} with $\eps=\eps_m$ for any $m\geq 1$, we have
that $G(\pi S^{-1} z_0)=0$.


\item[(C8)] If $G(\alpha)\neq 0$ for any
$\alpha\in \partial V$ and $d(G,V)\neq 0$, then there exists
$\eps_1>0$ sufficiently small such that for each
$\varepsilon\in(0,\eps_1]$  there is  at least one $T$-periodic
solution $\varphi_\varepsilon$ of system \eqref{ps} such that
$\rho(\varphi_\varepsilon(0),\mathcal{Z})\to 0$ as $\varepsilon\to
0$.
\end{itemize}
In addition we assume that there exists $\alpha_0\in V$ such that
$G(\alpha_0)=0$, $G(\alpha)\neq 0$ for all $\alpha\in \overline V
\setminus \{\alpha_0\}$ and $d(G,V)\neq 0$, and we denote
$z_0=S\left(
\begin{array}{c}\alpha_0\\\beta_0(\alpha_0)\end{array}\right)$.
Moreover we also assume:
\begin{itemize}
\item[(A8)] There exists ${\delta_1}>0$ and $L_G>0$ such that
\[ ||G(\alpha_1)-G(\alpha_2)||\geq L_G ||\alpha_1-\alpha_2||, \, \mbox{ for all } \alpha_1,\alpha_2\in B_{\delta_1}(\alpha_0),\]
\item[(A9)] For $\delta>0$ sufficiently small there exists
$M_\delta\subset [0,T]$ Lebesgue measurable with
$mes(M_\delta)=o(\delta)/\delta$  such that
\begin{equation*}
||g(t,z_1+\zeta,\eps)-g(t,z_1,0)-g(t,z_2+\zeta,\eps)+g(t,z_2,0)||\leq
o(\delta)/\delta ||z_1-z_2||\ ,
\end{equation*}
for all $t\in [0,T]\setminus M_\delta$ and for all $z_1,z_2\in
B_{\delta}(z_0)$, $\eps\in[0,\delta]$ and $\zeta\in B_\delta(0)$.
\end{itemize}
Then the following conclusion holds.
\begin{itemize}
\item[(C9)] There exists $\delta_2>0$  such that for any $\eps \in
(0,\eps_1]$, $\varphi_{\eps}$ is the only $T$-periodic solution of
\eqref{ps} with initial condition in $B_{\delta_2}(z_0)$. Moreover
$\varphi_{\eps}(0)\to z_0$ as $\eps\to 0$.
\end{itemize}
\end{thm}

\noindent To prove the theorem we need three preliminary lemmas that
are interesting by themselves. For example, in Lemma \ref{lem1} we
prove the existence of the derivative (in $\eps=0$) with respect to
some parameter denoted $\eps$ of the solution of some initial value
problem without assuming that the system is $C^1$. We also study the
properties of this derivative.


\begin{lem} \label{lemlip}
Let $f\in C^2(\mathbb{R}^n,\mathbb{R}^n)$ and $K_1,K_2$ be compact
subsets of $\mathbb{R}^n.$ Then the following inequality holds for
all $x^0_1,x^0_2\in K_1$, $y_1,y_2\in K_2$ and $\eps\in[0,1]$.
\begin{equation} \label{pty}
||f(x^0_1+\eps y_1)-f(x^0_1)-f(x^0_2+\eps y_2)+f(x^0_2)||\leq
O(\eps)||x^0_1-x^0_2||+O(\eps) ||y_1-y_2||\, .
\end{equation}
In addition  for $m>0$ sufficiently  small and $ u_1,\, u_2,\,
v_1,\, v_2\in B_m(0)\subset \mathbb{R}^n$ we have
\begin{equation} \label{ptu}
\begin{array}{ll}
\,||f\left( x^0_1+v_1 +\eps y^0_1+\eps u_1 \right) -
f\left( x^0_1+v_1 \right)- \eps f'(x^0_1)y^0_1-\\

\medskip

 f\left( x^0_2+ v_2+\eps y^0_2+\eps
u_2 \right)+f\left( x^0_2+ v_2 \right)+\eps f'(x^0_2)y^0_2|| \leq\\

\medskip

 \left[ o(\eps) + \eps O(m)
\right] ||x^0_1-x^0_2||+O(\eps) || v_1-v_2||+\\

\medskip

\left[ o(\eps)+\eps O(m)\right] ||
y^0_1-y^0_2||+O(\eps)||u_1-u_2||\, .
\end{array}
\end{equation}
\end{lem}

\noindent {\bf Proof.}  We define $\Phi(x^0,y,\eps)=f(x^0+\eps
y)-f(x^0)$ for all $x^0\in \overline{\mbox{co}}K_1$, $y\in
\overline{\mbox{co}}K_2$ and $\eps\in[0,1]$. Relation (\ref{pty})
follows from the mean value inequality applied to $\Phi_i$ with
$i\in\overline{1,n}$ and the following estimations.
\begin{eqnarray*}
\frac{\partial \Phi_i}{\partial x^0}(x^0,y,\eps)&=&
(f_i)^{\prime}(x^0+\eps y)-
 (f_i)^{\prime}(x^0)=O(\eps)\quad{\rm and}\\
\frac{\partial \Phi_i}{\partial y}(x^0,y,\eps)&=&\eps
(f_i)^{\prime}(x^0+\eps y)=O(\eps).
\end{eqnarray*}
In order to prove  relation (\ref{ptu}) we define
\[ \Phi(x^0,v,y^0,u,\eps)=f( x^0+v+\eps y^0+\eps u ) -
f( x^0+v)- \eps f'(x^0)y^0\, ,\] for all $x^0\in
\overline{\mbox{co}}K_1$, $y^0\in \overline{\mbox{co}}K_2$, $u,v\in
B_m(0)$ and $\eps\in[0,1]$. We apply again the mean value inequality
to the components $\Phi_i,$ $i\in\overline{1,n},$ using the
following estimations.
\begin{eqnarray*}
\frac{\partial \Phi_i}{\partial x^0}(x^0,v,y^0,u,\eps)&&= (f_i)^{\prime}( x^0+v+\eps y^0+\eps u )- (f_i)^{\prime}(x^0+v)- \eps (f_i)''(x^0)y^0\\
&&=o(\eps)+\eps (f_i)''(x^0+v)u\\&&+ \ \eps \left[ (f_i)''(x^0+v)-(f_i)''(x^0)\right] y^0\\&&=o(\eps)+\eps O(m)+\eps o(m)/m=o(\eps)+\eps O(m),\\
\frac{\partial \Phi}{\partial v}(x^0,v,y^0,u,\eps)&&=
(f_i)^{\prime}( x^0+v+\eps y^0+\eps u
)- (f_i)^{\prime}(x^0+v)=O(\eps),\\
\frac{\partial \Phi_i}{\partial y^0}(x^0,v,y^0,u,\eps)&&=\eps
(f_i)^{\prime}( x^0+v+\eps y^0+\eps u )- \eps
(f_i)^{\prime}(x^0)\\&&=\eps (f_i)^{\prime}( x^0+v+\eps y^0+\eps u
)- \eps (f_i)^{\prime}(x^0+v)\\
& & +\ \eps (f_i)^{\prime}(x^0+v)-\eps (f_i)^{\prime}(x^0)\\&&=o(\eps)+\eps O(m),\\
\frac{\partial \Phi}{\partial u}(x^0,v,y^0,u,\eps)&&=\eps
(f_i)^{\prime}( x^0+v+\eps y^0+\eps u )=O(\eps)\, .
\end{eqnarray*}
 $\Box$

\bigskip


\begin{lem}\label{lem1} We consider $f\in C^2(\mathbb{R}\times \mathbb{R}^n,\mathbb{R}^n)$
 and  $g\in C^0(\mathbb{R}\times \mathbb{R}^n\times[0,1],\mathbb{R}^n)$ a locally uniformly
 Lipschitz function
 with respect to the second variable. 
  For $z\in \mathbb{R}^n$ and $\eps\in [0,1]$, we denote by $\,x(\cdot,z,\eps)\,$  the unique solution of \[ \dot{x}=f(t,x)+\eps
  g(t,x,\eps)\, , \quad x(0)=z,\]  and by $\,y(t,z,\eps)= \left[x(t,z,\eps)-x(t,z,0)\right]/\eps\,$ $($here $\eps\neq 0$ $)$.
We assume that for a given $T>0$ there exist a compact set $K\subset
\mathbb{R}^n$ with nonempty interior and ${\delta}>0$ such that
$x(t,z,\eps)$ is well--defined for all $t\in [0,T]$, $z\in K$ and
$\eps\in[0,{\delta}]$.  Then the following statements hold.
\begin{itemize}
 \item[(C10)] There exists
$y(t,z,0)=\lim \limits_{\eps \to 0}y(t,z,\eps)$ being the solution
of the initial value problem
\begin{equation*} 
\dot{y}(t)=D_xf(t,x(t,z,0))y+g(t,x(t,z,0),0)\, ,\quad y(0)=0.
\end{equation*}
The above limit holds uniformly with respect to $(t,z)\in[0,T]\times
K$.
\item[(C11)] The functions $x,y:[0,T]\times K \times [0,\delta]\to
\mathbb{R}^n$  are continuous and uniformly Lipschitz with respect
to their second variable. \item[(C12)] In addition if there exists
$z_0\in
   \mbox{int}(K)$ such that assumption
{\em (A9)} of Theorem \ref{thm2} holds with the same small $\delta>0$ as above, then  
\begin{equation*} 
 ||y(t,z_1+\zeta,\eps)-y(t,z_1,0)-y(t,z_2+\zeta,\eps)+y(t,z_2,0)||\leq
o(\delta)/\delta ||z_1-z_2||\, ,
\end{equation*}
for all $t\in [0,T]$, $z_1,z_2\in B_\delta (z_0)$, $\eps \in
[0,\delta]$ and $\zeta \in B_\delta (0)$.
\end{itemize}
\end{lem}

\noindent {\bf Proof.} (C10) We define
$\displaystyle{\tilde{f}(t,z,\eps)=\frac{f(t,x(t,z,\eps))-f(t,x(t,z,0))}{x(t,z,\eps)-x(t,z,0)}}$
for $\eps\neq 0$ and $\tilde{f}(t,z,0)=D_xf(t,x(t,z,0))$. In this
way we obtain the continuous function $\tilde f:[0,T]\times K \times
[0,\delta]\to \mathbb{R}^n$. 
For $\eps\neq 0$, using the definitions of $x(t,z,\eps)$ and
$y(t,z,\eps)$ we deduce immediately that $y(0,z,\eps)=0$ and also
that \begin{equation} \label{eq:y} \dot{y}(t,z,\eps)=
\tilde{f}(t,x(t,z,\eps)) y(t,z,\eps) +
g(t,x(t,z,\eps),\eps).\end{equation} Passing to the limit as $\eps
\to 0$, we obtain that $y(\cdot,z,0)$ is the solution of the given
initial value problem. Hence \eqref{eq:y} holds also for $\eps=0$.
Since the right hand side of \eqref{eq:y} is given by a continuous
function, we have that the limit $y(t,z,0)=\lim \limits_{\eps \to
0}y(t,z,\eps)$ holds uniformly with respect to $(t,z)\in[0,T]\times
K$.

\medskip

(C11) The facts that the functions $x,y:[0,T]\times K \times
[0,\delta]\to \mathbb{R}^n$ are continuous, and that $x$ is
Lipschitz with respect to its second variable can be obtained as a
corollary of the general theorem on the dependence of the solutions
of an ordinary differential equation on the parameters (see
\cite[Lemma~8.2]{ama}).

 It remains to prove that $y:[0,T]\times K \times
[0,\delta]\to \mathbb{R}^n$ is uniformly Lipschitz with respect to
its second variable.

There exist compact subsets $K_1$ and $K_2$ of $\mathbb{R}^n$ such
that $x(t,z,\eps)\in K_1$ and $y(t,z,\eps)\in K_2$ for all
$(t,z,\eps)\in [0,T]\times K \times [0,\delta]$.

Moreover the representation $x(s,z,\eps)=x(s,z,0)+\eps y(s,z,\eps)$
allows to use Lemma \ref{lemlip}, relation (\ref{pty}) with
$x^0_1=x(s,z_1,0)$, $x^0_2=x(s,z_2,0),$ $y_1=y(s,z_1,\eps),$
$y_2=y(s,z_2,\eps)$ in order to obtain
\begin{eqnarray*}
 & \|f(t,x(t,z_1,\eps))-f(t,x(t,z_1,0))-f(t,x(t,z_2,\eps))+ f(t,x(t,z_2,0))\|\le \\
 &
 O(\eps)\|x(t,z_1,0)-x(t,z_2,0)\|+O(\eps)\|y(t,z_1,\eps)-y(t,z_2,\eps)\|,
\end{eqnarray*}
for all $t\in[0,T]$, $z\in K$ and $\eps\in[0,\delta]$.
 This last inequality and the fact that $g$ is locally uniformly
Lipschitz, used together with the representation
\begin{equation*}
 y(t,z,\eps)=\frac{1}{\eps}\int\limits_0^t\left[f(s,x(s,z,\eps))-f(s,x(s,z,0))\right]ds+\int\limits_0^t
 g(s,x(s,z,\eps),\eps)ds,\label{stst}
\end{equation*}
imply that
\begin{eqnarray*}
 \left\|y(t,z_1,\eps)-y(t,z_2,\eps)\right\|&\leq&
 O(\delta)/\delta
 \int\limits_0^t\left\|y(s,z_1,\eps)-y(s,z_2,\eps)\right\|ds \\ && +O(\delta)/\delta
 \int\limits_0^t\left\|x(s,z_1,\eps)-x(s,z_2,\eps)\right\|ds, \end{eqnarray*}
 for all $t\in [0,T]$, $z_1,z_2\in K$ and $\eps \in [0,\delta].$

\noindent We use now the fact that  the function
$
x(t,z,\varepsilon)$ is  Lipschitz with respect to  $z$ and we deduce
\[
  \left\|y(t,z_1,\eps)-y(t,z_2,\eps)\right\|\le
O(\delta)/\delta\left\|z_1-z_2\right\|+O(\delta)/\delta\int\limits_0^t\left\|y(s,z_1,\eps)-y(s,z_2,\eps)\right\|ds.
\]
Applying  Gr\"{o}nwall lemma (see \cite[Lemma~6.2]{hale} or
\cite[Ch.~2, Lemma \S~11]{dem}) we finally have for all $t\in
[0,T]$, $z_1,z_2\in K$, $\eps \in [0,\delta]$,
$||y(t,z_1,\eps)-y(t,z_2,\eps)||\leq
O(\delta)/\delta\left\|z_1-z_2\right\|$.

\medskip

(C12) First we note that assumption (A9) of Theorem 2 and the fact
that $g$ is locally uniformly Lipschitz with respect to the second
variable assure that the following relation holds
\begin{equation}
\label{gunic}
\begin{array}{ll}
||g(t,z_1+\zeta_1,\eps)-g(t,z_1,0)-g(t,z_2+\zeta_2,\eps)+g(t,z_2,0)||\leq\\
\leq o(\delta)/\delta||z_1-z_2||+O(\delta)/\delta ||\zeta_1-\zeta_2||,
\end{array}
\end{equation}
for all $t\in [0,T]\setminus M_{\delta}$, $z_1,z_2\in B_{\delta}
(z_0)$, $\eps \in [0,\delta]$ and $\zeta_1,\zeta_2 \in B_\delta
(0)$. We introduce the notations $ v(t,z,\zeta)=
x(t,z+\zeta,0)-x(t,z,0)$,
$\widetilde{\zeta}(s,z,\zeta,\eps)=v(s,z,\zeta)+\eps
y(s,z+\zeta,\eps)$ and
$u(t,z,\zeta,\eps)=y(t,z+\zeta,\eps)-y(t,z,0)$. Since the function
$x(\cdot,\cdot,0)$ is $C^1$, $v$ is Lipschitz with respect to $z$ on
$[0,T]\times K \times B_\delta(0)$ with some constant
$o(\delta)/\delta$, we have
\begin{eqnarray*}
&&u(t,z,\zeta,\eps)=y(t,z+\zeta,\eps)-y(t,z,0)\\&&=\frac 1{\eps}
\int _0^t \left[ f(s,x(s,z+\zeta,\eps))-f(s,x(s,z+\zeta,0)) - \eps
D_x f(s,x(s,z,0))y(s,z,0) \right] ds  \\ &&+ \int _0^t \left[
g(s,x(s,z+\zeta,\eps),\eps)-g(s,x(s,z,0),0) \right]ds.
\end{eqnarray*}
Our aim is to estimate a Lipschitz constant with respect to $z$ of the function $u$ on
 $[0,T]\times B_{\delta}(z_0)\times B_\delta(0)\times [0,\delta]$. We will apply Lemma \ref{lemlip}, relation
 (\ref{gunic}),
 the fact that $g$ is locally uniformly Lipschitz, and
 using  the following decompositions and estimations that hold for
 $(s,z,\zeta,\eps)\in [0,T]\times B_{\delta}(z_0)\times B_\delta(0)\times
 [0,\delta]$,
\begin{eqnarray*} x(s,z+\zeta,\eps)&=&x(s,z,0)+v(s,z,\zeta)+\eps y(s,z,0) +\eps
u(s,z,\zeta,\eps),\\
 x(s,z+\zeta,0)&=&x(s,z,0)+ v(s,z,\zeta),\\
x(s,z+\zeta,\eps)&=&x(s,z,0)+\widetilde{\zeta}(s,z,\zeta,\eps),\end{eqnarray*}
\[ ||v(t,z,\zeta)||\leq o(\delta)/\delta, \quad ||u(t,z,\zeta,\eps)||\leq o(\delta)/\delta, \quad
||\widetilde{\zeta}(s,z,\zeta,\eps)||\leq \delta O(\delta)/\delta,\]
we obtain
\begin{eqnarray*}
&&||u(t,z_1,\zeta,\eps) - u(t,z_2,\zeta,\eps)||\leq \\ && \frac 1
{\eps} \int_0^t  \left[ o(\eps)+\eps o(\delta)/\delta \right]
||x(s,z_1,0)-x(s,z_2,0)|| +O(\eps) ||v(s,z_1,\zeta)-v(s,z_2,\zeta)
||+ \\ && \left[ o(\eps)+\eps o(\delta)/\delta \right]
||y(s,z_1,0)-y(s,z_2,0)||+O(\eps)|| u(s,z_1,\zeta,\eps) -
u(s,z_2,\zeta,\eps)|| ds+\\ && \int_{(0,t)\setminus M_{\delta}}
o(\delta)/\delta||x(s,z_1,0)-x(s,z_2,0)||
+O(\delta)/\delta||\widetilde{\zeta}(s,z_1,\zeta,\eps)-\widetilde{\zeta}
(s,z_2,\zeta,\eps)||ds+\\&& o(\delta)/\delta ||z_1-z_2|| .
\end{eqnarray*}
Now we use that some Lipschitz constants with respect to $z$ for the
functions $x$ and $y$ on  $[0,T]\times B_{\delta}(z_0)\times
[0,\delta]$ are $O(\delta)/\delta$, while for the functions $v$ on
$[0,T]\times B_{\delta}(z_0)\times [0,\delta]$ and
$\widetilde{\zeta}$ on $[0,T]\times B_{\delta}(z_0)\times
B_\delta(0)\times [0,\delta]$ are $o(\delta)/\delta$, and finally we
obtain that
\begin{eqnarray*}
&&||u(t,z_1,\zeta,\eps) - u(t,z_2,\zeta,\eps)||\leq\\&&
o(\delta)/\delta ||z_1-z_2|| + O(\delta)/\delta \int _0^t
||u(t,z_1,\zeta,\eps) - u(t,z_2,\zeta,\eps)||ds\,.
\end{eqnarray*}
The conclusion follows after applying the Gr\"{o}nwall inequality.

\qed

\begin{lem}\label{lem2}
We consider the $C^1$ function $Y$ acting from $\mathbb{R}^n$ into
the space of $n\times n$ matrices,  the $C^2$ function
  $\widetilde{P}: \mathbb{R}^n\to\mathbb{R}^n$ and $z_*\in \mathbb{R}^n$ such that
$\widetilde{P}(z_*)=0$. We denote $P:\mathbb{R}^n\to\mathbb{R}^n$
the $C^1$ function given by
 $P(z)=Y(z)\widetilde{P}(z)$ for all $z\in \mathbb{R}^n$. Then
$DP(z_*)=Y(z_*)D\widetilde{P}(z_*)$, $P$ is twice differentiable in
$z_*$ and, for each $i\in \overline{1,n}$, the Hessian matrix
$HP_i(z_*)$ is symmetric.
\end{lem}

\noindent {\bf Proof.}  We
have $\displaystyle{DP(z)=\left( \frac{\partial Y}{\partial z_1}(z)\widetilde{P}(z),...,
\frac{\partial Y}{\partial z_n}(z)\widetilde{P}(z)\right)+Y(z)D\widetilde{P}(z)}$ for all
$z\in \mathbb{R}^n$. From this it follows the formula for $DP(z_*)$ since $\widetilde{P}(z_*)=0$.\\

In order to prove that  $P$ is twice differentiable in $z_*$, taking
into account the above expression of $DP$, it is enough to prove
that for each
 $i\in \overline{1,n}$, $\displaystyle{z\mapsto \frac{\partial Y}{\partial
z_i}(z)\widetilde{P}(z)}$  and $z\mapsto Y(z)D\widetilde{P}(z)$ are
differentiable in $z_*$. The last map is $C^1$, hence
 it remains to prove the differentiability only for the first one. We fix  $i\in \overline{1,n}$.
From the relation
$$
\begin{array}{l}
\displaystyle{ \frac{\partial Y}{\partial z_i}(z_*+h)\widetilde{P}(z_*+h)-\frac{\partial Y}{\partial z_i}(z_*)\widetilde{P}(z_*)}=\\
\displaystyle{ \frac{\partial Y}{\partial
z_i}(z_*+h)\left(\widetilde{P}(z_*+h)-\widetilde{P}(z_*)\right)
 = \frac{\partial Y}{\partial z_i}(z_*+h)D\widetilde{P}(z_*)+o(h)},
 \end{array}
$$
we deduce that $\displaystyle{z\mapsto \frac{\partial Y}{\partial
z_i}(z)\widetilde{P}(z)}$ is differentiable in $z_*$ and that\\
$\displaystyle{D\left(\frac{\partial Y}{\partial z_i}\cdot
\widetilde{P}\right)(z_*)=\frac{\partial Y}{\partial
z_i}(z_*)D\widetilde{P}(z_*).}$

 In order to prove that the Hessian matrix $HP_i(z_*)$ is symmetric,
 for every $j,\, k\in \{1,...,n\}$
we must prove that
\[ \frac{\partial^2 P_i}{\partial z_j \partial z_k}(z_*)=\frac{\partial^2 P_i}{\partial z_k \partial
z_j}(z_*).\] We denote by $Y_i(z)$ the $i$--th row of the $n\times
n$ matrix $Y(z)$. For all $z\in \mathbb{R}^n$ we have
\[ \frac{\partial P_i}{\partial z_j}(z)=Y_{i}(z)\frac{\partial \widetilde{P}}{\partial
z_j}(z)+\frac{\partial Y_{i}}{\partial
z_j}(z)\widetilde{P}(z).\] Then
\begin{equation*}
\frac{\partial^2 P_i}{\partial z_j \partial z_k}(z_*)=
\frac{\partial Y_{i}}{\partial z_k}(z_*)\frac{\partial
\widetilde{P}}{\partial z_j}(z_*)+Y_{i}(z_*)\frac{\partial
^2\widetilde{P}}{\partial z_j\partial z_k}(z_*)+\frac{\partial
Y_{i}}{\partial z_j}(z_*)\frac{\partial \widetilde{P}}{\partial
z_k}(z_*)\,.
\end{equation*}
Since $\widetilde{P}$ is $C^2$ it is easy to check the symmetry of
this last relation with respect to $(j,k)$. \qed\\

\noindent {\bf Proof of Theorem~\ref{thm2}.} We need to study the
zeros of the function $z\mapsto x(T,z,\varepsilon)-z\ ,$ or
equivalently of
$$F(z,\varepsilon)=Y^{-1}(T,z)(x(T,z,\varepsilon)-z).$$
 The function $F$ is well defined at least for any $z$ in some small neighborhood of
 $\mathcal{Z}$ and any $\eps\geq 0$ sufficiently small.  We will apply Theorem 1. We denote
\[P(z)=Y^{-1}(T,z)\left( x(T,z,0)-z \right),\quad Q(z,\eps)=Y^{-1}(T,z)y(T,z,\eps),
\]
where  $y(t,z,\eps)=[x(t,z,\eps)-x(t,z,0)]/\eps$, like in Lemma
\ref{lem1}. Hence $F(z,\eps)=P(z)+\eps Q(z,\eps)$.

The fact that $f$ is $C^2$ assures that  the function $z\mapsto
x(T,z,0)$ is also $C^2$ (see \cite[Ch.~4, \S~24]{pont}).
 Since (see \cite[Ch.~III, Lemma
\S~12]{dem}) $\left(Y^{-1}(\cdot,z)\right)^*$ is a fundamental
matrix solution of the system
$$
  \dot u=-(D_x f(t,x(t,z,0),0))^*u\, ,
$$
and $f$ is $C^2$,  we have that the matrix function $(t,z)\mapsto
\left(Y^{-1}(t,z)\right)^*$ is  $C^1$. Therefore the matrix function
$(t,z)\mapsto Y^{-1}(t,z)$, and consequently also the function $P$
are  $C^1$.

By Lemma~\ref{lem1} we now conclude that $Q$ is continuous, locally
uniformly Lipschitz with respect to $z$, and
\begin{equation} \label{eq:qz0} Q(z,0)=\int _0^T Y^{-1}(s,z)g(s,x(s,z,0),0) ds.\end{equation}
Since, by our hypothesis (A6), $x(\cdot,z,0)$ is $T$-periodic for
all $z\in\mathcal{Z}$ we have that $x(T,z,0)-z=0$ for all $z\in
\mathcal{Z}$, and consequently $P\left( z\right)=0$ for all $z\in
\mathcal{Z}$. This means that hypothesis (A1) of Theorem 1 holds.
Moreover applying Lemma \ref{lem2} we have that
\begin{equation*}
  DP(z)=Y^{-1}(T,z)\left(\frac{\partial x}{\partial
  z}(T,z,0)-I_{n\times n}\right)\quad{\rm for \ any\
}z\in\mathcal{Z},
\end{equation*}
and $P$ satisfies hypothesis (A3) of Theorem \ref{thm1}.
 But  $\displaystyle{\left({\partial x}/{\partial
z}\right)(\cdot,z,0)}$ is the normalized fundamental matrix of the
linearized system (\ref{ls}) (see \cite[Theorem 2.1]{kraope}).
Therefore $\displaystyle{\left({\partial x}/{\partial
z}\right)(t,z,0)}=Y(t,z)Y^{-1}(0,z)$, and we can write
\begin{equation}\label{zg}
DP\left( z\right)=Y^{-1}(0,z)-Y^{-1}(T,z)\quad{\rm for \ any\
}z\in\mathcal{Z}.
\end{equation}
Using our hypothesis (A7) we see that also assumption (A2) of
Theorem~\ref{thm1} is satisfied. From the definition of $G$ and
relation (\ref{eq:qz0}) we have that
$$
  G(\alpha)=\pi Q\left(S\left(\begin{array}{c}\alpha \\ \beta_0(\alpha)\end{array}
  \right),0\right)\,.
$$
That is, the function denoted in  Theorem~\ref{thm1} by
$\widehat{Q}$ is here $G$, and it satisfies the hypotheses of
Theorem~\ref{thm1}. Moreover, note that when $G$ satisfies (A8) then
assumption (A4) of Theorem \ref{thm1} is fulfilled.

\medskip

(C7) Follows from (C1) of Theorem~\ref{thm1}.

\medskip

(C8)  Follows from (C2) of Theorem~\ref{thm1}.


\medskip

(C9) In order to prove the uniqueness of the $T$-periodic solution,
it remains only to check (A5) of Theorem \ref{thm1}. For doing this
we show that the function
 $(z,\zeta,\eps)\in B_{\delta}(z_0)\times B_\delta(0)\times [0,\delta]\mapsto
Q(z+\zeta,\eps)-Q(z,0)$  is Lipschitz with respect to $z$ with
some constant $o(\delta)/\delta$.  We write
\begin{eqnarray*}
Q(z+\zeta,\eps)-Q(z,0)&=&Y^{-1}(T,z+\zeta)\left[ y(T,z+\zeta
,\eps)-y(T,z,0)\right]+\\&&\left[
Y^{-1}(T,z+\zeta)-Y^{-1}(T,z)\right]y(T,z,0)\,.
\end{eqnarray*}
It is known that for proving that a sum of two functions is
Lipschitz with some constant $o(\delta)/\delta$, it is enough  to
prove that each function is Lipschitz with such constant; while in
order to prove that a product of two functions is Lipschitz with
some constant $o(\delta)/\delta$, it is sufficient to prove that
both functions are Lipschitz and only one of them is bounded by some
constant $o(\delta)/\delta$ and Lipschitz with respect to $z$ with
some constant $o(\delta)/\delta$.

By Lemma \ref{lem1} we know  that the function $z\in
B_{\delta}(z_0)\mapsto y(T,z,0)$ is Lipschitz. The fact that
$z\mapsto Y^{-1}(T,z)$ is $C^1$ assures that $(z,\zeta)\in
B_{\delta}(z_0)\times B_\delta(0)\mapsto Y^{-1}(T,z+\zeta)$ is
Lipschitz with respect to $z$.

From Lemma \ref{lem1} we have that the function
\[(z,\zeta,\eps)\in B_\delta(z_0)\times B_\delta(0)\times
[0,\delta]\mapsto  y(T,z+\zeta ,\eps)-y(T,z,0)\] is bounded by some
constant $o(\delta)/\delta$ and Lipschitz with some constant
$o(\delta)/\delta$. Since $z\mapsto Y^{-1}(T,z)$ is $C^1$, the same
is true for the function \[(z,\zeta)\in B_{\delta}(z_0)\times
B_\delta(0)\mapsto \left[ Y^{-1}(T,z+\zeta)-Y^{-1}(T,z)\right].\]
Hence $Q$ satisfies (A5) of Theorem \ref{thm1} and the conclusion
holds.

\qed\\

By using Theorem~\ref{thm2} we can provide now a result which
includes both the existence-uniqueness results by Malkin
\cite{mal} and by Melnikov \cite{mel}. The main contribution of
our result is that we do not impose any assumptions on the
algebraic multiplicity of the multiplier $+1$ of (\ref{nnn}).
Also, the condition imposed on the function $z\mapsto g(t,z,\eps)$
is weaker than the $C^1$ condition used by Malkin and Melnikov. In
particular our Theorem~\ref{thm3} covers the class of piecewise
differentiable systems.

\begin{thm} \label{thm3} Assume that $f\in C^2(\mathbb{R}\times\mathbb{R}^n, \mathbb{R}^n)$ and
 $g\in C^0(\mathbb{R}\times\mathbb{R}^n\times[0,1],\mathbb{R}^n)$ are $T$-periodic in the first
variable, and  that $g$ is locally uniformly Lipschitz with respect
to the second variable. 
Assume that the unperturbed system \eqref{np} satisfies the
following conditions.
\begin{itemize}
\item[(A10)]
There exists an open ball $U\subset \mathbb{R}^k$ with $k\leq n$ and
a function $\xi\in C^1(\overline U,\mathbb{R}^n)$ such that for any
$h\in \overline U$ the $n\times k$ matrix $D\xi(h)$ has rank $k$ and
$\xi(h)$ is the initial condition of a $T$-periodic solution of
\eqref{np}.
***\item[(A11)] For each $h\in \overline U$ the linear system \eqref{ls} with $z=\xi(h)$  has the Floquet multiplier $+1$ with the geometric
multiplicity
 equal to $k$.
\end{itemize}
 Let $u_1(\cdot,h), \dots ,u_k(\cdot,h)$ be linearly independent $T$-periodic solutions of the
 adjoint linear
system
\begin{equation}\label{ss}
  \dot u=-\left(D_x f(t,x(t,\xi(h),0))\right)^*u,
\end{equation}
such that $u_1(0,h), \dots ,u_k(0,h)$ are $C^1$ with respect to  $h$
*** and define the function $M:\overline U \to \mathbb{R}^k$ (called
the Malkin's bifurcation function) by
$$
  M(h)=\int\limits_0^T \left(\begin{array}{c}
    \left<u_1(s,h),g(s,x(s,\xi(h),0),0)\right>\\ ... \\
    \left<u_k(s,h),g(s,x(s,\xi(h),0),0)\right>\end{array}\right)
    ds.
$$
Then the following statements hold.
\begin{itemize}
\item[(C13)] For any sequences $(\varphi_m)_{m\geq 1}$ from $C^0(\mathbb{R}, \mathbb{R}^n)$ and $(\eps_m)_{m\geq 1}$
from $[0,1]$ such that $\varphi_m(0)\to \xi(h_0)\in \xi(\overline
U)$, $\eps_m\to 0$ as $m\to \infty$ and $\varphi_m$ is a
$T$-periodic solution of \eqref{ps} with $\eps=\eps_m$, we have that
$M(h_0)=0$.

\item[(C14)] If $M(h)\neq 0$ for any
$h\in \partial U$ and $d(M,U)\neq 0$, then there exists $\eps_1>0$
 sufficiently small such that for each $\eps\in (0,\eps_1]$ there is at least one $T$-periodic solution $\varphi_\varepsilon$ of
 system \eqref{ps}  such that
$\rho(\varphi_\varepsilon(0),\xi(\overline U))\to 0$ as
$\varepsilon\to 0$.
\end{itemize}

In addition we assume that there exists $h_0\in U$ such that
$M(h_0)=0$, $M(h)\neq 0$ for all $h\in \overline U \setminus
\{h_0\}$ and $d(M,U)\neq 0$. Moreover we  assume that hypothesis
{\em (A9)} of Theorem \ref{thm2} holds with $z_0=\xi(h_0)$ and that
\begin{itemize}
\item[(A12)] There exists $\delta_1>0$ and  $L_M>0$ such that
\[ ||M(h_1)-M(h_2)||\geq L_M ||h_1-h_2||, \mbox{ for all } h_1,h_2\in B_{\delta_1}(h_0).\]
\end{itemize}
Then the following conclusion holds.
\begin{itemize}
\item[(C15)] There exists $\delta_2>0$ such that for any $\eps\in (0,\eps_1]$, $\varphi_\eps$ is the only $T$-periodic solution of
\eqref{ps} with initial condition in $B_{\delta_2}(z_0)$. Moreover
$\varphi _{\eps}(0)\to \xi(h_0)$ as $\eps\to 0$.
\end{itemize}
\end{thm}

***\noindent {\bf Remark.} {\it The existence of $k$ linearly
independent $T$-periodic solutions of the adjoint linear system
\eqref{ss} follows by hypothesis {\em (A10)} (see e.g.
\cite[Ch.~III, \S~23, Theorem~2]{dem}). Indeed, we have that
$y_i(t,h)=D_zx(t,\xi(h),0)D_{h_i}\xi(h)$
 for  $i\in\overline{1,k}$ are solutions of \eqref{ls} and they are
 linearly independent on the base of {\em (A10)}. The assertion follows by
 the fact that a linear system and its adjoint have the same number
 of linearly independent solutions. Moreover, hypothesis {\em
 (A11)} assures that there is no other $T$-periodic solution to
 \eqref{ls} linearly independent of these. } ***\\

***\noindent {\bf Proof.} We apply Theorem \ref{thm2}. For the moment
we describe the set $\mathcal{Z}$ that appear in hypothesis (A6) as
$\mathcal{Z}=\bigcup\limits_{h\in\overline{U}}
\left\{\xi(h)\right\}.$ First we find the matrix $S$  such that
hypothesis (A7) holds. In order to achieve this, for each
$z\in\mathcal{Z}$ we denote $U(t,z)$ some fundamental matrix
solution of \eqref{ss} that has in its first $k$ columns the
$T$-periodic solutions $u_1, \dots ,u_k$ and such that $U(0,z)$ is
$C^1$. Then the first $k$ columns of the matrix $U(0,z)-U(T,z)$ are
null vectors. The matrix $Y(t,z)$ such that $Y^{-1}(t,z)=[U(t,z)]^*$
is a fundamental matrix solution of \eqref{ls}, i.e. of the system
($z=\xi(h)\in\mathcal{Z}$)
\begin{equation}\label{lsh}
  \dot y=D_x f(t,x(t,\xi(h),0))y.
\end{equation} Then the first $k$ lines of the matrix
$Y(0,z)^{-1}-Y(T,z)^{-1}$ are null vectors. Since the Floquet
multiplier $1$ of \eqref{ls} has geometric multiplicity $k$ we have
that the matrix $Y^{-1}(0,z)-Y^{-1}(T,z)$ has range $n-k$. Hence
this matrix has $n-k$ linearly independent columns. We claim that
there exists an invertible matrix $S$ such that the matrix
$\left(Y^{-1}(0,z)-Y^{-1}(T,z)\right)S$ has in the first $k$ lines
null vectors and in the lower right corner some $(n-k)\times (n-k)$
 invertible matrix $\Delta(z)$.  With this we prove that (A7) holds.
 In order to justify the claim we note first that whatever the matrix $S$ would
 be, the first $k$ lines of $\left(Y^{-1}(0,z)-Y^{-1}(T,z)\right)S$
 are null vectors. Now we choose an invertible matrix $S$
 such that its last $(n-k)$ columns are vectors of the form
$$
  e_i=\left(\begin{array}{c} 0_{(i-1)\times 1}\\ 1\\
  0_{(n-i)\times 1}\end{array}\right),\quad i\in\overline{1,n}
$$
distributed in such a way that the $n-k$ linearly independent
columns of $Y^{-1}(0,z)-Y^{-1}(T,z)$ become the last $n-k$ columns
of $\left(Y^{-1}(0,z)-Y^{-1}(T,z)\right)S$. Now it is easy to see
that the $(n-k)\times (n-k)$ matrix from the lower right corner of
$\left(Y^{-1}(0,z)-Y^{-1}(T,z)\right)S$ is invertible.

Now we come back to prove (A6). By taking the derivative with
respect to $h\in U$ of $ \dot{x}(t,\xi(h))=f(t,x(t,\xi(h)))$
 we obtain that $D_\xi x(\cdot,\xi(h))\cdot D\xi(h)$ is a matrix solution for (\ref{lsh}).
 But $x(\cdot,\xi(h))$ is $T$-periodic for any
$h\in U,$ therefore $D_\xi x(\cdot,\xi(h))\cdot D\xi(h)$ is
$T$-periodic. This fact assures that each column of $D\xi(h)$ is the
initial condition of some $T$-periodic solution of (\ref{lsh}) and
these $T$ periodic solutions are the columns of
$Y(t,\xi(h))Y^{-1}(0,\xi(h))D\xi(h)$. Then
$Y(T,\xi(h))Y^{-1}(0,\xi(h))D\xi(h)=D\xi(h)$, that further gives
$\left[ Y^{-1}(0,\xi(h)) -Y^{-1}(T,\xi(h))\right]S S^{-1}D\xi(h)=0$.
Hence the columns of $S^{-1}D\xi(h)$ belong to the kernel of $\left[
Y^{-1}(0,\xi(h)) -Y^{-1}(T,\xi(h))\right]S$. Since  (A7) holds we
have that the kernel of $\left[ Y^{-1}(0,\xi(h))
-Y^{-1}(T,\xi(h))\right]S$ contains vectors whose last $n-k$
components are null. We deduce that  there exists some $k\times k$
matrix, denoted by $\Psi$, such that
\begin{equation}\label{nu1}
 S^{-1}D\xi(h)=\left(\begin{array}{c} \Psi \\ 0_{(n-k)\times k}
  \end{array}\right).
\end{equation}
Since by the assumption (A10) the matrix $D\xi(h)$ has rank $k$
 and $S^{-1}$ is
invertible, we have that the matrix $S^{-1}D\xi(h)$ 
 should also have rank $k$,
  that is only possible if
\begin{equation}\label{nu2}
  {\rm det} \Psi \not=0.
\end{equation}
 We fix some $h_*\in U$ and we denote
$\alpha_*=\pi S^{-1}\xi(h_*).$ Using (\ref{nu1}) and (\ref{nu2}),
and applying the
    Implicit Function Theorem we have that there exists an open ball, neighborhood of $\alpha_*$, denoted
    ${V}\subset\mathbb{R}^k$, and a $C^1$ function
    $\widetilde{h}:\overline V\to {U}$  such that
    \begin{equation} \label{eq:implicit}
      \pi S^{-1}\xi (\widetilde{h}(\alpha))=\alpha\quad {\rm for\ any} \quad \alpha\in{\overline V}.
    \end{equation}
Now we define the $C^1$ function $\beta_0:{\overline V}\to
\mathbb{R}^{n-k}$ as $
  \beta_0(\alpha)=\pi^\bot S^{-1}\xi(\widetilde{h}(\alpha)).
$ Note that $S\left(\begin{array}{c}\alpha \\
\beta_0(\alpha)\end{array}\right)=\xi(\widetilde{h}(\alpha))$.
 Hence the assumption (A6) of Theorem~\ref{thm2} is satisfied with
$S$, $V$ and $\beta_0$ defined as above.

The bifurcation function $G$ defined in Theorem~\ref{thm2} can be
written using our notations  as
$$
  G(\alpha)=\pi\int\limits_0^T
  Y^{-1}(s,\xi(\widetilde{h}(\alpha)))g(s,x(s,\xi(\widetilde{h}(\alpha)),0))ds.
$$
Since $Y^{-1}(s,\xi(h))=[U(t,h)]^*$ (see the beginning of the proof)
we have that in the first $k$ lines of $Y^{-1}(s,\xi(h))$ are the
vectors $(u_1(s,h))^*$, ... , $(u_n(s,h))^*$ and so
$$
  G(\alpha)=\int\limits_0^T \left(\begin{array}{c}
    \left<u_1(s,{\widetilde{h}(\alpha)}),\,g(s,x(s,\xi(\widetilde{h}(\alpha)),0),0)\right>\\ ... \\
    \left<u_k(s,{\widetilde{h}(\alpha)}),\,g(s,x(s,\xi(\widetilde{h}(\alpha)),0),0)\right>\end{array}\right)
    ds.
$$
From here one can see that there is the following relation between $G$ and the Malkin bifurcation function $M$,
\begin{equation} \label{eq:G} G(\alpha)=M\left( \tilde{h}(\alpha)\right) \mbox{ for any } \alpha\in \overline V.\end{equation}
***

(C13) Follows from (C7) of Theorem~2.

(C14) Without loss of generality (we can diminish $U,$ if necessary)
we can consider  that $\widetilde{h}$ is a homeomorphism from $V$
onto $U$ and taking into account that $C$ is an invertible matrix by
\cite[Theorem 26.4]{krazab} we have
$$ {\rm deg}(G,V)={\rm deg}(M,U).$$ Thus (C14) follows  applying
conclusion (C8) of Theorem~\ref{thm2}.\\

(C15) We need only to prove  assumption (A8) of Theorem~\ref{thm2}
provided that our hypothesis (A12) holds.
 First, taking the derivative of (\ref{eq:implicit}) with respect to $\alpha$ and using (\ref{nu1}), we obtain that
$D\widetilde{h}(\alpha_*)=\Psi^{-1}$, hence it is invertible and,
moreover, $L_h=\|D\widetilde{h}(\alpha_*)\|/2\neq 0$. We have that
there exists $\delta>0$ sufficiently small such that
$\|D\widetilde{h}(\alpha)-D\widetilde{h}(\alpha_*)\|\leq L_h$ for
all $\alpha\in B_\delta(\alpha_*)$. Using the generalized Mean Value
Theorem (see \cite[Proposition 2.6.5]{clark}), we have that 
\[||\widetilde{h}(\alpha_1)-\widetilde{h}(\alpha_2)||\geq L_h
||\alpha_1-\alpha_2||\,  \mbox{ for all } \alpha_1,\alpha_2\in
B_\delta(\alpha_0).\] Since $C$ is invertible, $M$ satisfies (A10)
and (\ref{eq:G}), we deduce that $G$ satisfies hypothesis (A8) of
Theorem 2.

\qed

\section{An example}
\label{sec4} In this section we illustrate Theorem~\ref{thm3} by
studying the existence of $2\pi$-periodic solutions for the
following four--dimensional nonsmooth
 system
\begin{equation} \label{system}
\begin{array}{lll}
\dot{x_1}&=&\,\,\,\,\,x_2-x_1(x_1^2+x_2^2-1),\\
\dot{x_2}&=&-x_1-x_2(x_1^2+x_2^2-1)+\eps \left(\sin t + \varphi (k_1x_3)\right),\\
\dot{x_3}&=&\,\,\,\,\,x_4-x_3(x_3^2+x_4^2-1),\\
\dot{x_4}&=&-x_3-x_4(x_3^2+x_4^2-1)+\eps  \varphi (k_2x_2),
\end{array}
\end{equation}
when $\eps>0$ is sufficiently small, $k_1$ and $k_2$ are arbitrary reals and
$\varphi:\mathbb{R}\to \mathbb{R}$ is the piecewise linear function
\begin{equation} \nonumber
\varphi(x)=\left\{
\begin{array}{ll}
-1,~~~\mbox{ for }x\in (-\infty,-1),\\
\,\,\,\,\,x,~~~~\mbox{ for }x\in [-1,1],\\
\,\,\,\,\,1,~~~~\mbox{ for }x\in (1,\infty).
\end{array}
\right.
\end{equation}

Before proceeding to this study we write down a sufficient condition
in order that the function $g$ satisfies (A9). This is of general
interest, not only for this example. It was also used in \cite{blm}.
\begin{itemize}
\item[(A13)] For any $\delta>0$ sufficiently small there exists
$M_\delta\subset [0,T]$ Lebesgue measurable with mes$(M_\delta)=o(\delta)/\delta$
and such that for every $t\in [0,T]\setminus M_\delta$ and for all $z\in
B_\delta(z_0)$, $\eps\in[0,\delta]$ ,
$ \,
||D_zg(t,z,\eps)-D_zg(t,z_0,0)||\leq
o(\delta)/\delta \,. $
\end{itemize}
The fact that (A13) implies (A9) follows from the Mean Value Theorem.\\

 Coming back to system (\ref{system}) we define $g:\mathbb{R}\times \mathbb{R}^4\to \mathbb{R}^4$ as
  $g(t,z_1,z_2,z_3,z_4)=(0,\,\sin t + \varphi (k_1z_3),\,0, \,\varphi (k_2z_2))$. This function is uniformly locally Lipschitz and satisfies (A13) (hence also (A9)) for any $z_0=(z_1^0,z_2^0,z_3^0,z_4^0)$ with $|k_2z_2^0|\neq 1$ and $|k_1z_3^0|\neq 1$.\\

We study now the unperturbed system. For $\eps=0$  system (\ref{system}) becomes
\begin{eqnarray*}
\dot{x_1}&=&x_2-x_1(x_1^2+x_2^2-1), ~~~~~~\,\,\,\,\,\,~~~~~~\dot{x_2}=-x_1-x_2(x_1^2+x_2^2-1),\\
\dot{x_3}&=&x_4-x_3(x_3^2+x_4^2-1),
~~~~~~~\,\,\,\,\,\,~~~~~\dot{x_4}=-x_3-x_4(x_3^2+x_4^2-1),
\end{eqnarray*}
and it has a family of $2\pi$-periodic orbits, whose initial
conditions are given by
\[\xi(\theta,\eta)=(\sin \theta, \cos \theta, \sin \eta, \cos \eta),\quad (\theta,\eta)\in \mathbb{R}^2.\]
We have $\,x(t,\xi(\theta,\eta),0)=\left( \sin (t+\theta),
\cos(t+\theta), \sin (t+\eta), \cos(t+\eta)\right).$

It is easy to see that the function $\xi:\mathbb{R}^2\to
\mathbb{R}^4$ satisfies assumption (A10) of Theorem~\ref{thm3}. We
consider now the linearized system,
\begin{equation} \label{lin}
\dot{y}=D_xf(t,x(t,\xi(\theta,\eta),0))y,
\end{equation}
where
$f(x_1,x_2,x_3,x_4)=(x_2-x_1(x_1^2+x_2^2-1),-x_1-x_2(x_1^2+x_2^2-1),x_4-x_3(x_3^2+x_4^2-1),-x_3-x_4(x_3^2+x_4^2-1))$,
and the following $2\pi$-periodic matrix
\[ \Phi(t,(\theta,\eta))=\left(
\begin{array}{ll}
-\cos(t+\theta) \hspace{1cm} 0 \hspace{1cm} \sin(t+\theta) \hspace{1cm} 0\\
\,\,\,\,\,\,\sin(t+\theta)
 \hspace{1cm} 0 \hspace{1cm} \cos(t+\theta) \hspace{1cm} 0\\
\hspace{1cm}0 \hspace{0.9cm} -\cos(t+\eta)       \hspace{1cm} 0           \hspace{1cm}\sin(t+\eta)\\
\hspace{1cm}0\hspace{1.4cm} \sin(t+\eta)       \hspace{1cm}0
\hspace{0.5cm}-\cos(t+\eta)
\end{array}
\right).\] Denoting
$\hat{\Phi}(t,(\theta,\eta))=\Phi(t,(\theta,\eta))\Phi^{-1}(0,(\theta,\eta))$,
it can be checked that the normalized fundamental matrix of system
(\ref{lin}) is
\[ \widehat{Y}(t,\xi(\theta,\eta))=\hat{\Phi}(t,(\theta,\eta))\left(
\begin{array}{ll}
I_{2\times 2} \hspace{0.8cm} 0_{2\times2}\\
\hspace{0.3cm}0_{2\times2} \hspace{0.5cm} e^{-2t}I_{2\times 2}
\end{array}
\right),\] and that it satisfies the hypothesis (A11) of
Theorem~\ref{thm3}.

A couple of linearly independent $2\pi$-periodic solutions of the
adjoint system to system (\ref{lin}) is
\begin{eqnarray*}
u_1(t)&=&\left( -\cos (t+\theta),\,\,\sin (t+\theta),\,\,0,\,\,0\right)^T,\\
u_2(t)&=&\left( 0,\,\,0,\,\,-\cos (t+\eta),\,\,\sin (t+\eta)\right)^T.
\end{eqnarray*}
Therefore the Malkin's bifurcation function $M(\theta,\eta)$ takes the form
$$ M(\theta,\eta)=
\left(\begin{array}{c} M_1(\theta,\eta)\\
M_2(\theta,\eta)\end{array}\right)=\left(\begin{array}{c}
\int\limits_0^{2\pi} \sin (s+\theta)\left[ \sin s+ \varphi\left(
k_1\sin
(s+\eta)\right)\right]ds\\
\int\limits_0^{2\pi} \sin (s+\eta) \varphi\left( k_2\cos
(s+\theta)\right)ds\end{array}\right).
$$
We denote
\begin{eqnarray*}
I(k)=\int_0^{2\pi} \sin t \, \varphi(k\sin t)dt, \quad i_1=I(k_1), \\ J(k)=\int_0^{2\pi} \cos t \, \varphi(k\cos
t)dt,  \quad  j_2=J(k_2).
\end{eqnarray*}
Then
\begin{eqnarray*}
M_1(\theta,\eta)&=& \pi \cos \theta + i_1\cos(\theta-\eta),\\
M_2(\theta,\eta)&=&  -j_2\sin (\theta-\eta).
\end{eqnarray*}
It is easy to see that the function $M$ is differentiable and that
the determinant of its Jacobian at each point $(\theta,\eta)$ is
\[\det\left(DM(\theta,\eta)\right)=-\pi j_2\sin \theta \cos(\theta-\eta).\]
Studying the zeros of the function $M$, we obtain the following results.

\medskip

If $j_2\neq 0$ and $|i_1|< \pi$, $M$ has exactly $4$ zeros, namely \[(a_1,\pi+a_1),\quad  (\pi+a_1,a_1),
\quad (\pi-a_1,\pi-a_1), \quad (2\pi-a_1,2\pi-a_1),\]
 and all have nonvanishing index (here
$a_1=\arccos (i_1/\pi)$).

If $j_2=0$, we have that $M_2\equiv 0$, hence
 the zeros are not isolated.

If $j_2\neq 0$ and $|i_1|=\pi$, $M$ has exactly $2$ zeros, namely
\[(a_1,\pi+a_1), \quad (\pi-a_1,\pi-a_1), \]
 and both have index $0$ (note that $a_1$ is $0$ or $\pi$).

  If  $j_2\neq 0$ and $|i_1|>\pi$, $M$ has no zeros.

\medskip

In order to complete the study, we calculate
\[ I(k)=J(k)=\left\{
\begin{array}{ll}
k\pi,~~~~0\leq k \leq 1,\\
2k\arcsin \frac 1k +2 \frac{\sqrt{k^2-1}}{k}, ~~~~k>1.
\end{array}
\right.\] We note that $I$ and $J$ are odd functions.  For $k>1$,
$I^{\prime}(k)=\arcsin \frac 1k -\frac 1k \sqrt{1-\frac
1{k^2}}>\frac 1k -\frac 1k \sqrt{1-\frac 1{k^2}}>0$. Then
$I(k)>I(1)=\pi$ for all $k>1$. In fact $|I(k)|>\pi$ if and only if
$|k|>1$. Also  $|I(k)|<\pi$ if and only if $|k|<1$, $|I(k)|=\pi$ if
and only if $|k|=1$, and $J(k)=0$ if and only if
$k=0$.\\

One can see that when $k_2\neq 0$ and $|k_1|<1$ the bifurcation function $M$ has exactly four zeros. The values of the function $\xi$ on these zeros are of the form $(\pm \sin a_1, \, \pm \cos a_1,\,\pm \sin a_1, \, \pm \cos a_1)$ where $a_1=\arccos k_1$.\\

Using all these facts and applying Theorem~\ref{thm3}, we obtain the
following result concerning the existence of $2\pi$--solutions of
system (\ref{system}) when $k_2\neq 0$.

\begin{pro} Let $k_2\neq 0$ and $\eps>0$ be sufficiently small.

 If $|k_1|>1$  system \eqref{system} has no $2\pi$-periodic
solutions with initial conditions converging to some point of
$\xi\left( \mathbb{R}\times \mathbb{R}\right)$ as $\eps\to 0$.

 If
$|k_1|<1$  system \eqref{system} has at least four $2\pi$-periodic
solutions with initial conditions converging to exactly four points
of $\xi\left( \mathbb{R}\times \mathbb{R}\right)$ as $\eps\to 0$. If
$|k_1k_2|\neq 1$ then the obtained $2\pi$-periodic solutions are
isolated.
\end{pro}

For $k_2=0$ the behavior of system (\ref{system}) is the same as of
the uncoupled system whose last two equations remain the same and
the first two change to
\begin{equation} \label{system2}
\begin{array}{lll}
\dot{x_1}&=&\,\,\,\,\,x_2-x_1(x_1^2+x_2^2-1),\\
\dot{x_2}&=&-x_1-x_2(x_1^2+x_2^2-1)+\eps \left(\sin t + \varphi
(k_1\sin(t+\eta))\right).
\end{array}
\end{equation}
In the following we study the existence of $2\pi$-periodic solutions
of the above planar system. We consider only the case that
$(k_1,\eta) \in \mathbb{R}\times \mathbb{R}\setminus \{ (1,\pi),
(-1,0)\}$, otherwise the perturbation term is identically zero. Note
also that the perturbation term does not depend on the space
variable. System (\ref{system2}) with $\eps=0$ has a family of
$2\pi$-periodic orbits whose initial conditions are given by
\[\xi(\theta)=(\sin \theta,\cos\theta),\quad \theta\in \mathbb{R}.\]
For each fixed $\eta$ the Malkin's bifurcation function is
\[M^\eta(\theta)=\pi \cos \theta + i_1\cos(\theta-\eta),\] and it
can be easily seen that it has exactly $2$ zeros, and both have
nonvanishing index. If $k_1=0$ or $\eta\in \{ 0,\pi\}$, the zeros of
$M^\eta$ are $\pi/2$ and $3\pi/2$. Otherwise the zeros are
$\pi-\theta^*$ and $2\pi-\theta^*$, where
$\displaystyle{\theta^*=\arctan \frac{\pi+i_1\cos
\eta}{i_1\sin\eta}.}$

 Hence applying Theorem \ref{thm3} we obtain the result.
\begin{pro} Let $(k_1,\eta) \in \mathbb{R}\times \mathbb{R}\setminus \{ (1,\pi),
(-1,0)\}$ and $\eps>0$ be sufficiently small.  Then system
\eqref{system2} has  two isolated  $2\pi$-periodic solutions with
initial conditions converging to exactly two points of
$\xi^\eta\left( \mathbb{R}\right)$ as $\eps\to 0$.
\end{pro}

\end{document}